# The Relationship Between Two Commutators

Keith A. Kearnes     Ágnes Szendrei [*]

**Abstract**

We clarify the relationship between the linear commutator and the ordinary commutator by showing that in any variety satisfying a nontrivial idempotent Mal'cev condition the linear commutator is definable in terms of the centralizer relation. We derive from this that abelian algebras are quasi–affine in such varieties. We refine this by showing that if **A** is an abelian algebra and $\mathcal{V}(\mathbf{A})$ satisfies an idempotent Mal'cev condition which fails to hold in the variety of semilattices, then **A** is affine.

## 1 Introduction

In general algebra, the last twenty years have been marked by the development and application of various general commutator theories. Each theory defines a binary operation, called a "commutator", on congruence lattices of certain kinds of algebras. The operation might be considered to be a generalization of the usual commutator operation of group theory. The first development of a general commutator theory is due to Jonathan Smith in [19]. His theory is essentially a complete extension of the group commutator to any variety of algebras which has permuting congruences. This theory was then extended to congruence modular varieties by Hagemann and Herrmann [7]. Alternate ways of extending the theory to congruence modular varieties are described by Gumm [6] and Freese and McKenzie [3].

Some new definitions of the congruence modular commutator were suggested in the 1980's and early 1990's, and it was discovered that some of these definitions made sense, although different kinds of sense, for classes of objects which are not necessarily congruence modular and which are not necessarily varieties and where the objects are not necessarily algebras. Some examples of this can be found in [2], [9], [12], [14], [16], [18], [24]. The most successful of these theories, measured in terms of the number of significant theorems proved which do not mention the commutator operation itself, is clearly the commutator theory developed in [9] based on the term condition. This commutator is often called "the TC commutator" and we shall adhere to this convention in this introduction. In subsequent sections of this paper it will simply be referred to as *the* commutator.

The virtues of the TC commutator are that it is the easiest to calculate with and that it is effective in describing the connection between the way operations of an algebra compose and the distribution of its congruences. The shortcomings of the TC commutator are that it has poor categorical properties and that centrality is not understood well. In particular, we are far from understanding the structure of abelian algebras.

---

[*]Research supported by the Hungarian National Foundation for Scientific Research grant no. T 17005

An attempt to remedy the first shortcoming of the TC commutator (poor categorical properties) can be found in [2]. Here the commutator of two congruences is defined to be their free intersection. This approach has not had much success yet, but at least the authors of [2] are able to prove that the free intersection "is a true commutator", where we take this phrase to mean that it "agrees with the usual commutator in any congruence modular variety".

An attempt to remedy the second of the shortcomings of the TC commutator is hinted at in [18]. In some sense this paper *begins* by stating a structure theorem for abelian algebras and then working backwards to find the commutator that gives this structure theorem! More precisely, [18] starts off by noticing that the leading candidate for what an abelian algebra *should be* is an algebra which is representable as a subalgebra of a reduct of an affine algebra. Then a characterization of such algebras is given. The paper [18] hints that there is a commutator associated with this approach, and that commutator has come to be known as the "linear commutator".

For the last ten years the outstanding open question concerning the linear commutator has been the following: Is the linear commutator "a true commutator"? In other words, does the linear commutator agree with the usual commutator in every congruence modular variety? (This is Problem 5.18 of [13].) We answer this question affirmatively in this paper. We show moreover that the linear commutator agrees with the commutator defined in [3], which is a symmetric commutator operation defined by the term condition, in any variety which satisfies a nontrivial idempotent Mal'cev condition. Our result is the sharpest possible connection between the linear commutator and the term condition for varieties which satisfy a nontrivial idempotent Mal'cev condition. The coincidence of the linear commutator with a commutator defined by the term condition promises that we can have both a good understanding of centrality and the ease of calculating with the term condition in such varieties.

The most important immediate corollary of our result on the linear commutator is that it shows that an abelian algebra is quasi–affine if it generates a variety satisfying a nontrivial idempotent Mal'cev condition. We extend this to show that an abelian algebra is affine if it generates a variety satisfying an idempotent Mal'cev condition which fails in the variety of semilattices. This fact is a significant extension of Herrmann's Theorem (see [8]), which is itself a cornerstone of modular commutator theory. Further corollaries of our main theorem include the following: (i) Abelian algebras are affine in any variety which satisfies a nontrivial lattice identity as a congruence equation. (This is an affirmative answer to the question asked in Problem 1 of [9].) (ii) The property of having a weak difference term is equivalent to a Mal'cev condition. (iii) The congruence equation $[\alpha, \beta] = \alpha \wedge \beta$ is equivalent to congruence meet semidistributivity.

## 2  Several Commutators

In this section we will define the commutator operations which play a role in this paper. All the definitions are based on properties of the "$\alpha, \beta$–matrices", which we now define.

**Definition 2.1** Let $\mathbf{A}$ be an algebra and assume that $\alpha, \beta \in \mathrm{Con}(\mathbf{A})$. $M(\alpha, \beta)$ is the set of all matrices
$$\begin{bmatrix} t(\mathbf{a}, \mathbf{b}) & t(\mathbf{a}, \mathbf{b}') \\ t(\mathbf{a}', \mathbf{b}) & t(\mathbf{a}', \mathbf{b}') \end{bmatrix}$$
where $t(\mathbf{x}, \mathbf{y})$ be an $(m + n)$–ary term in the language of $\mathbf{A}$, $\mathbf{a}, \mathbf{a}' \in A^m$, $\mathbf{b}, \mathbf{b}' \in A^n$, and $(a_i, a_i') \in \alpha$, $(b_i, b_i') \in \beta$.

We say that $\alpha$ **centralizes** $\beta$ **modulo** $\delta$ and write $C(\alpha, \beta; \delta)$ if whenever
$$\begin{bmatrix} a & b \\ c & d \end{bmatrix} \in M(\alpha, \beta),$$
then $a \equiv_\delta b$ implies $c \equiv_\delta d$. To connect this with what was said in the introduction, we say that the $\alpha, \beta$–**term condition** holds if $C(\alpha, \beta; 0)$ holds. It is not hard to see that, for a fixed $\alpha$ and $\beta$, the set of all $\delta$ such that $C(\alpha, \beta; \delta)$ holds is closed under complete intersection. Therefore, there is a least $\delta$ such that $C(\alpha, \beta; \delta)$ holds. The class of all $\delta$ such that both $C(\alpha, \beta; \delta)$ and $C(\beta, \alpha; \delta)$ hold is also closed under complete intersection, so there is a least $\delta$ in this set, too.

**Definition 2.2** Let $\mathbf{A}$ be an algebra and assume that $\alpha, \beta \in \mathrm{Con}(\mathbf{A})$. The **commutator** of $\alpha$ and $\beta$, written $[\alpha, \beta]$, is the least $\delta$ such that $C(\alpha, \beta; \delta)$ holds. The **symmetric commutator** of $\alpha$ and $\beta$, written $[\alpha, \beta]_s$, is the least $\delta$ such that $C(\alpha, \beta; \delta)$ and $C(\beta, \alpha; \delta)$ hold.

Obviously we have $[\alpha, \beta] \leq [\alpha, \beta]_s = [\beta, \alpha]_s$. Moreover, it is easy to see that both operations are monotone in each variable.

In [3], Freese and McKenzie develop a commutator theory for congruence modular varieties which is based on the symmetric commutator. In their Proposition 4.2 they prove that the equation $[\alpha, \beta] = [\alpha, \beta]_s$ holds in every congruence modular variety. Clearly this equation holds for an algebra exactly when the operation $[\alpha, \beta]$ is symmetric. From this it is easy to see that the equation $[\alpha, \beta] = [\alpha, \beta]_s$ holds in any variety with a difference term (see [10]), but that it fails already in some varieties with a weak difference term, such as the variety of inverse semigroups. The point of this remark is just to note that the two commutator operations we have defined can be different, but they agree in "nice" varieties.

Our next task is to define the linear commutator. For any similarity type $\tau$ we let $\tau^*$ denote the expansion of $\tau$ to include the new symbol $p$. If $\mathcal{V}$ is the variety of all $\tau$–algebras, then we let $\mathcal{V}^*$ be the variety of all algebras of type $\tau^*$ such that

(1) $p(x, y, y) = x = p(y, y, x)$,

(2) $p$ commutes with itself, and

(3) for any basic $\tau$–operation $f$ and for tuples of variables $\mathbf{u}$ and $\mathbf{v}$ we have
$$f(\mathbf{u}, p(x, y, z), \mathbf{v}) = p(f(\mathbf{u}, x, \mathbf{v}), f(\mathbf{u}, y, \mathbf{v}), f(\mathbf{u}, z, \mathbf{v})).$$

Conditions (1) and (2) imply that for any $\mathbf{A} \in \mathcal{V}^*$ we have $p(x, y, z) = x - y + z$ with respect to some abelian group structure on $A$. We call a structure of the form $\langle A; p \rangle$ where $p$ satisfies (1) and (2) an **affine abelian group**. Condition (3) says that every basic $\tau$–operation is multilinear with respect to $p$. This does not imply that all $\tau$–terms are multilinear. For example, if $\tau$ has a single binary operation $xy$, then the term $b(x, y) = (xy)x$ is linear in $y$ but it is not linear in $x$ for some members of $\mathcal{V}^*$. What can be proved by induction on the complexity of a term $t(\mathbf{x})$ is that if the parse tree of $t$ has the property that no two leaves are labelled with the same variable, then $t(\mathbf{x})$ is multilinear. Therefore, every $\tau$–term can be obtained from a multilinear term by identifying variables. For example, $b(x, y)$ can be obtained from the multilinear term $(xy)z$ by identifying $z$ with $x$.

Note that Conditions (1), (2) and (3) are equational and so $\mathcal{V}^*$ is a variety which, by (1), is congruence permutable.

There is a forgetful functor $F : \mathcal{V}^* \to \mathcal{V}$ which "forgets $p$", and this functor has a left adjoint $* : \mathcal{V} \to \mathcal{V}^*$. It is not difficult to describe this adjoint. For any $\mathbf{A} \in \mathcal{V}$ one can take for the universe of $\mathbf{A}^*$ the free affine abelian group generated by $A$, for $p$ one can take the group operation $x - y + z$, and for each basic $\tau$–operation one can take the multilinear extension of the corresponding operation of $\mathbf{A}$. The morphism part of the functor $*$ can be described as follows: If $\varphi : \mathbf{A} \to \mathbf{B}$ is a $\mathcal{V}$–homomorphism, then $\varphi^* : \mathbf{A}^* \to \mathbf{B}^*$ is simply the linear extension of $\varphi$. We leave it to the reader to verify that $\mathbf{A}^* \in \mathcal{V}^*$ and that restriction to $\mathbf{A}$ is a natural bijection from $\mathrm{Hom}_{\mathcal{V}^*}(\mathbf{A}^*, \mathbf{B})$ to $\mathrm{Hom}_{\mathcal{V}}(\mathbf{A}, F(\mathbf{B}))$. This verifies that $\langle *, F \rangle$ is an adjoint pair.

If $\alpha \in \mathrm{Con}(\mathbf{A})$, then the natural homomorphism $\varphi : \mathbf{A} \to \mathbf{A}/\alpha$ extends to a homomorphism $\varphi^* : \mathbf{A}^* \to (\mathbf{A}/\alpha)^*$ and the range of $\varphi^*$ contains $A/\alpha$, which is a generating set for $(\mathbf{A}/\alpha)^*$, therefore $\varphi^*$ is surjective. We let $\alpha^*$ denote the kernel of $\varphi^*$. By the First Isomorphism Theorem we have $\mathbf{A}^*/\alpha^* \cong (\mathbf{A}/\alpha)^*$.

**Definition 2.3** Let $\mathbf{A}$ be an algebra and assume that $\alpha, \beta \in \mathrm{Con}(\mathbf{A})$. The **linear commutator** of $\alpha$ and $\beta$, written $[\alpha, \beta]_\ell$, is $[\alpha^*, \beta^*]|_A$.

Since $\mathbf{A}^*$ belongs to the congruence permutable variety $\mathcal{V}^*$, $[\alpha^*, \beta^*] = [\alpha^*, \beta^*]_s$, and therefore the linear commutator could as easily be defined by $[\alpha, \beta]_\ell = [\alpha^*, \beta^*]_s|_A$.

The fact that the commutator is symmetric in a congruence permutable variety implies that the linear commutator is also symmetric. The fact that the commutator is monotone in each variable and that the mapping $\alpha \mapsto \alpha^*$ is monotone implies that the linear commutator is monotone in each variable.

**LEMMA 2.4** *Let $\mathbf{A}$ be an algebra. The following hold.*

(1) $\alpha^*|_A = \alpha$.

(2) $[\alpha, \beta]_s \leq [\alpha, \beta]_\ell \leq \alpha \wedge \beta$.

(3) *The pair of mappings $\langle f, g \rangle$ where $f : \alpha \mapsto \alpha^*$ and $g : \beta \mapsto \beta|_A$ is an adjunction from $\mathrm{\mathbf{Con}}(\mathbf{A})$ to $\mathrm{\mathbf{Con}}(\mathbf{A}^*)$.*

**Proof:** Item (1) follows from the fact that $\varphi^*|_A = \varphi$.

For the first part of (2) let $\delta = [\alpha^*, \beta^*] \in \mathrm{Con}(\mathbf{A}^*)$. Since $\mathbf{A}$ is a subalgebra of a reduct of $\mathbf{A}^*$, the fact that $C(\alpha^*, \beta^*; \delta)$ holds implies that $C(\alpha^*|_A, \beta^*|_A; \delta|_A)$ holds. From part (1) we get that $C(\alpha, \beta; [\alpha, \beta]_\ell)$ holds. The linear commutator is symmetric, so interchanging the roles of $\alpha$ and $\beta$ we get that $C(\beta, \alpha; [\alpha, \beta]_\ell)$ holds too. It follows that $[\alpha, \beta]_s \leq [\alpha, \beta]_\ell$. For the second part of (2) simply restrict the inequality $[\alpha^*, \beta^*] \leq \alpha^* \wedge \beta^*$ to $\mathbf{A}$.

There are two ways to interpret statement (3) depending on how you consider an ordered set to be a category; we specify the category by saying that there is a (unique) homomorphism from $\gamma$ to $\delta$ if $\gamma \leq \delta$. In this sense, what we must prove is that for any $\alpha \in \mathrm{Con}(\mathbf{A})$ and $\beta \in \mathrm{Con}(\mathbf{A}^*)$ we have
$$\alpha^* \leq \beta \iff \alpha \leq \beta|_A.$$
The forward implication holds since if $\alpha^* \leq \beta$, then by (1) we have $\alpha = \alpha^*|_A \leq \beta|_A$. If the reverse implication fails, then it fails if we replace $\beta$ by the intersection of all $\gamma$ for which $\alpha \leq \gamma|_A$. Therefore, we only need to show that if $\beta$ is the least congruence on $\mathbf{A}^*$ for which $\alpha \leq \beta|_A$, then $\alpha^* \leq \beta$. Equivalently, we must show that $\alpha^*$ is the least congruence on $\mathbf{A}^*$ whose restriction to $\mathbf{A}$ is $\alpha$. We show slightly more. By definition, $\alpha^*$ is the kernel of the natural surjective homomorphism $\mathbf{A}^* \to (\mathbf{A}/\alpha)^*$. But since $\mathbf{A}^*$ and $(\mathbf{A}/\alpha)^*$ are freely generated as affine abelian groups by $A$ and $A/\alpha$ respectively, and since this natural homomorphism respects the affine abelian group operation $p$, it follows that $\alpha^*$ is none other than the affine abelian group congruence on $\mathbf{A}^*$ generated by $\alpha$. Thus $\alpha^*$ is the least equivalence relation on $\mathbf{A}^*$ which is compatible with $p(x, y, z)$ and which contains $\alpha$. A fortiori we have that $\alpha^*$ is the least congruence on $\mathbf{A}^*$ whose restriction to $\mathbf{A}$ contains $\alpha$. Thus (3) is proved. $\square$

If $[\alpha, \beta]_s \neq [\alpha, \beta]_\ell$, then by Lemma 2.4 (2) we have $[\alpha, \beta]_s < [\alpha, \beta]_\ell$. By Lemma 2.4 (1) and (3) we have $\gamma^* < [\alpha^*, \beta^*]$, where $\gamma = [\alpha, \beta]_s$. Therefore in $\mathbf{A}^*/\gamma^*$ we have $0 < [\alpha^*/\gamma^*, \beta^*/\gamma^*]$. But in the isomorphism between $\mathbf{A}^*/\gamma^*$ and $(\mathbf{A}/\gamma)^*$ we have that $\alpha^*/\gamma^*$ corresponds to $(\alpha/\gamma)^*$. This can be seen by considering the kernels of the natural surjections and their composite in:
$$\mathbf{A}^* \longrightarrow (\mathbf{A}/\gamma)^* \longrightarrow (\mathbf{A}/\alpha)^*.$$
Thus we get $0 < [(\alpha/\gamma)^*, (\beta/\gamma)^*]$ in $(\mathbf{A}/\gamma)^*$, implying that $0 < [\alpha/\gamma, \beta/\gamma]_\ell$ in $\mathbf{A}/\gamma$. At the same time we have $0 = \gamma/\gamma = [\alpha/\gamma, \beta/\gamma]_s$ in $\mathbf{A}/\gamma$. This proves the following corollary.

**COROLLARY 2.5** *If $\mathbf{A}$ is an algebra in which the symmetric commutator is different than the linear commutator, then a homomorphic image of $\mathbf{A}$ has a pair of congruences $\alpha$ and $\beta$ such that*
$$[\alpha, \beta]_s = 0 < [\alpha, \beta]_\ell.$$

Our next task is to characterize the situation where $[\alpha, \beta]_\ell = 0$. We model our arguments on pages 323–324 of [18].

If $\alpha$ is a congruence on the algebra $\mathbf{A}$, we write $\mathbf{A} \times_\alpha \mathbf{A}$ to denote the subalgebra of $\mathbf{A} \times \mathbf{A}$ which consists of those pairs whose first coordinate is $\alpha$–related to the second. That is, $\mathbf{A} \times_\alpha \mathbf{A}$ is the subalgebra with universe $\alpha$. If $\beta$ is also a congruence on $\mathbf{A}$, then we consider $M(\alpha, \beta)$ to be the binary relation on $\mathbf{A} \times_\alpha \mathbf{A}$ defined by

$$\begin{bmatrix} a \\ c \end{bmatrix} \text{ is related to } \begin{bmatrix} b \\ d \end{bmatrix} \iff \begin{bmatrix} a & b \\ c & d \end{bmatrix} \in M(\alpha, \beta).$$

This relation is clearly a tolerance on $\mathbf{A} \times_\alpha \mathbf{A}$. In the case where our algebra is $\mathbf{A}^*$ and our congruences are of the form $\alpha^*$ and $\beta^*$, we have that $M(\alpha^*, \beta^*)$ is transitive (since $\mathbf{A}^*$ has a Mal'cev operation), and therefore it is a congruence. Combining Lemma 4.8 of [3] with Theorem 4.9 (ii) of [3], we get that $(u,v) \in [\alpha^*, \beta^*]$ if and only if

$$\begin{bmatrix} u & u \\ u & v \end{bmatrix} \in M(\alpha^*, \beta^*).$$

Therefore, $[\alpha, \beta]_\ell > 0$ if and only if there are distinct $u, v \in A$ which determine a matrix in $M(\alpha^*, \beta^*)$ in this way. In order to write this only in terms of the algebra $\mathbf{A}$ we need to understand $M(\alpha^*, \beta^*)$ in terms of $\mathbf{A}$.

Fix an element of $A$ and label it 0. We will write $x + y$ to mean $p(x, 0, y)$, $-x$ to mean $p(0, x, 0)$ and $x - y$ to mean $x + (-y)$. These operations are abelian group operations on $A$ and with regard to these definitions we have $p(x, y, z) = x - y + z$

**LEMMA 2.6** (1) The elements of $\mathbf{A}^*$ are expressible as sums of elements of $A$, $\sum n_i a_i$, where $\sum n_i = 1$.

(2) Elements $w$ and $z$ of $\mathbf{A}^*$ are $\alpha^*$–related if and only if $w - z = \sum(a^j - b^j)$ where $(a^j, b^j) \in \alpha$.

(3) A $\tau^*$–term $F(\mathbf{x})$ is $\mathcal{V}^*$–equivalent to a sum $\sum n_i f_i(\mathbf{x})$ where each $f_i$ is a $\tau$–term and $\sum n_i = 1$.

(4) A matrix $M$ belongs to $M(\alpha^*, \beta^*)$ if and only if it equals a sum $\sum n_i N_i$ where each $N_i \in M(\alpha, \beta)$ and $\sum n_i = 1$.

(5) $(u, v) \in [\alpha, \beta]_\ell$ if and only if $u, v \in A$ and there exist matrices

$$\begin{bmatrix} a_i & b_i \\ c_i & d_i \end{bmatrix} \in M(\alpha, \beta)$$

such that $v - u = \sum(a_i - b_i - c_i + d_i)$.

**Proof:** Statement (1) follows from the fact that the universe of $\mathbf{A}^*$ is the free affine abelian group generated by the set $A$.

Statement (2) is a consequence of the fact (proved in our argument for Lemma 2.4 (3)) that $\alpha^*$ is the least equivalence relation on $\mathbf{A}^*$ which is compatible with $p(x, y, z)$ and which contains $\alpha$.

Statement (3) is proved by induction on the complexity of $F$ using the multilinearity of each basic $\tau$–operation. (We mention that although $f$ is a $\tau$–term, each component of $\mathbf{x}$ in $f(\mathbf{x})$ is a variable which ranges over elements of $\mathbf{A}^*$, not $\mathbf{A}$, and moreover some of these variables might be "fictitious".)

For (4), a typical matrix $M \in M(\alpha^*, \beta^*)$ is of the form

$$M = \begin{bmatrix} F(\mathbf{u}, \mathbf{w}) & F(\mathbf{u}, \mathbf{z}) \\ F(\mathbf{v}, \mathbf{w}) & F(\mathbf{v}, \mathbf{z}) \end{bmatrix}$$

where $F$ is a $\tau^*$–term, $(u_i, v_i) \in \alpha^*$ and $(w_i, z_i) \in \beta^*$. To prove (4) it suffices, by (3), to consider only the case where $F = f$ is a $\tau$–term. Furthermore, since every $\tau$–term is the specialization of a multilinear $\tau$–term, there is no loss of generality if we assume that $f$ is multilinear. Using (2), we can write $u_i = v_i + \sum(a_i^j - b_i^j)$ and $w_i = z_i + \sum(c_i^j - d_i^j)$ where $(a_i^j, b_i^j) \in \alpha$ and $(c_i^j, d_i^j) \in \beta$. Furthermore, by (1) we have that $v_i = \sum m_j^i g_j^i$ and $z_i = \sum n_j^i h_j^i$ where $\sum m_j^i = \sum n_j^i = 1$ and $g_j^i, h_j^i \in A$. Using the linearity of $f$ in its first variable we can expand $M = \begin{bmatrix} f(u_1-, -) & f(u_1-, -) \\ f(v_1-, -) & f(v_1-, -) \end{bmatrix}$ as

$$\begin{bmatrix} f(v_1-, -) & f(v_1-, -) \\ f(v_1-, -) & f(v_1-, -) \end{bmatrix} + \sum \left( \begin{bmatrix} f(a_1^j-, -) & f(a_1^j-, -) \\ f(a_1^j-, -) & f(a_1^j-, -) \end{bmatrix} - \begin{bmatrix} f(b_1^j-, -) & f(b_1^j-, -) \\ f(a_1^j-, -) & f(a_1^j-, -) \end{bmatrix} \right).$$

We can further expand the matrix on the left by substituting $\sum n_j^1 g_j^1$ for $v_1$. This has the effect of replacing the matrix on the left in the next line with the sum on the right.

$$\begin{bmatrix} f(v_1-, -) & f(v_1-, -) \\ f(v_1-, -) & f(v_1-, -) \end{bmatrix} = \sum n_j^1 \begin{bmatrix} f(g_j^1-, -) & f(g_j^1-, -) \\ f(g_j^1-, -) & f(g_j^1-, -) \end{bmatrix}.$$

Thus we may reduce $M$ to a sum of matrices which belong to $M(\alpha^*, \beta^*)$, each involving the $\tau$–term $f$, but where the entries in the first variable all belong to $A$ rather than to $\mathbf{A}^*$. Furthermore the sum of the coefficients of these matrices is 1. We can continue this process variable–by–variable replacing each matrix in this sum by a sum of matrices. At each step a single matrix is replaced by a sum of matrices in $M(\alpha^*, \beta^*)$ whose coefficients sum to 1. Therefore the sum of coefficients does not change during this process, it is always 1. Finally, when we are finished we have a sum of matrices of the form

$$N = \begin{bmatrix} f(\mathbf{p}, \mathbf{r}) & f(\mathbf{p}, \mathbf{s}) \\ f(\mathbf{q}, \mathbf{r}) & f(\mathbf{q}, \mathbf{s}) \end{bmatrix}$$

where $f$ is a $\tau$–term, $(p_i, q_i) \in \alpha^*|_A = \alpha$ and $(r_i, s_i) \in \beta^*|_A = \beta$. Such matrices belong to $M(\alpha, \beta)$. This proves one direction of (4). The other direction is easy since $M(\alpha^*, \beta^*)$ contains $M(\alpha, \beta)$ and is closed under $p(x, y, z) = x - y + z$ (since it is a subalgebra of $(\mathbf{A}^*)^4$). Thus any sum of the form $\sum n_i N_i$ with $\sum n_i = 1$ and $N_i \in M(\alpha, \beta)$ belongs to $M(\alpha^*, \beta^*)$.

To prove (5), assume that $(u, v) \in [\alpha, \beta]_\ell$. Then $(u, v) \in [\alpha^*, \beta^*]$, so

$$\begin{bmatrix} u & u \\ u & v \end{bmatrix} \in M(\alpha^*, \beta^*).$$

By (4) this means that

$$\begin{bmatrix} u & u \\ u & v \end{bmatrix} = \sum n_i N_i = \sum n_i \begin{bmatrix} a_i & b_i \\ c_i & d_i \end{bmatrix}$$

where all matrices shown belong to $M(\alpha, \beta)$. By replacing each $n_i N_i$ with $|n_i|$ copies of $N_i$ we may assume that each $n_i = \pm 1$. This leads to four coordinate equations: $u = \sum n_i a_i$, $u = \sum n_i b_i$, $u = \sum n_i c_i$ and $v = \sum n_i d_i$. Therefore we get

$$v - u = u - u - u + v = \left(\sum n_i a_i\right) - \left(\sum n_i b_i\right) - \left(\sum n_i c_i\right) + \left(\sum n_i d_i\right)$$

which may be written as $v - u = \sum n_i(a_i - b_i - c_i + d_i)$, where each $n_i = \pm 1$. Now observe that if some $n_i = -1$, then since

$$\begin{bmatrix} a'_i & b'_i \\ c'_i & d'_i \end{bmatrix} := \begin{bmatrix} b_i & a_i \\ d_i & c_i \end{bmatrix} \in M(\alpha, \beta)$$

and $(a'_i - b'_i - c'_i + d'_i) = -(a_i - b_i - c_i + d_i) = n_i(a_i - b_i - c_i + d_i)$, it is clear that we can reduce to the case where all $n_i = +1$. That is, we can find matrices in $\begin{bmatrix} a_i & b_i \\ c_i & d_i \end{bmatrix} \in M(\alpha, \beta)$ such that $v - u = \sum(a_i - b_i - c_i + d_i)$, proving one direction of (5).

For the other direction of (5), assume that $u, v \in A$, $\begin{bmatrix} a_i & b_i \\ c_i & d_i \end{bmatrix} \in M(\alpha, \beta)$ for each $i$ and that $v - u = \sum(a_i - b_i - c_i + d_i)$. Then the matrix

$$\begin{bmatrix} u & u \\ u & v \end{bmatrix} = \begin{bmatrix} u & u \\ u & u \end{bmatrix} + \sum\left(\begin{bmatrix} a_i & a_i \\ a_i & a_i \end{bmatrix} - \begin{bmatrix} a_i & b_i \\ a_i & b_i \end{bmatrix} - \begin{bmatrix} a_i & a_i \\ c_i & c_i \end{bmatrix} + \begin{bmatrix} a_i & b_i \\ c_i & d_i \end{bmatrix}\right)$$

is in $M(\alpha^*, \beta^*)$ by statement (4). Therefore $(u, v) \in [\alpha, \beta]_\ell$. This proves (5). □

Now we explain how to use Lemma 2.6 (5) to determine from $\mathbf{A}$ if $[\alpha, \beta]_\ell = 0$. This equality is equivalent to the implication $(u, v) \in [\alpha, \beta]_\ell \Longrightarrow u = v$, which may be written as: if $u, v \in A$, each $\begin{bmatrix} a_i & b_i \\ c_i & d_i \end{bmatrix} \in M(\alpha, \beta)$ and $v - u = \sum(a_i - b_i - c_i + d_i)$, then $u = v$. Since $u, v, a_i, b_i, c_i, d_i \in A$, the equality $v - u = \sum(a_i - b_i - c_i + d_i)$ holds in $\mathbf{A}^*$, the free affine abelian group with basis $A$, if and only if it holds in a trivial fashion: each element of $A$ in the expression $\sum(a_i - b_i - c_i + d_i)$ which occurs with a plus sign is matched with an identical element with a minus sign, with two elements left over which are equal to $-u$ and $+v$. Therefore, the condition that $[\alpha, \beta]_\ell = 0$ may be rephrased as follows.

> Given any (finite) sum of the form $\sum(a_i - b_i - c_i + d_i)$, where each summand comes from an $\alpha, \beta$–matrix, and where each element with a plus sign is matched with an identical element with a minus sign except that two elements $-u$ and $+v$ are left over, it is the case that $u = v$.

We can formulate this condition in terms of matchings between "positive elements" and "negative elements".

We find it notationally more convenient in the arguments that follow if we deal with $2 \times 2$ matrices whose first row contains both of the the positive elements and whose second row contains the negative elements, so let $TM(\alpha, \beta)$ denote the set of all matrices $\begin{bmatrix} a & d \\ c & b \end{bmatrix}$ such that $\begin{bmatrix} a & b \\ c & d \end{bmatrix} \in M(\alpha, \beta)$. (We will refer to $TM(\alpha, \beta)$ as the set of "twisted $\alpha, \beta$–matrices".)
In a twisted $\alpha, \beta$–matrix the first row contains what we have called the positive elements (of the sum $a - b - c + d$) and the second row contains the negative elements. Two elements in the same column are $\alpha$–related and we imagine a directed edge going *upward* from a bottom element to the top element in the same column. Such an edge will be called an $\alpha$–**edge**. Two elements on a diagonal are $\beta$–related and we imagine a directed edge going upward from a bottom element to a top element on the same diagonal. Such an edge will be called a $\beta$–**edge**.

Each twisted $\alpha, \beta$–matrix may be viewed as a copy of the graph $G$ depicted in Figure 1, with the vertices labelled with elements of $A$. We do not consider the vertex labelling to be part of the definition of $G$. Note that the labels that occur on the two vertices that determine an $\alpha$–edge must be $\alpha$–related and that the labels that occur on the two vertices that determine a $\beta$–edge must be $\beta$–related. However not every such labelling arises as a twisted $\alpha, \beta$–matrix, usually.

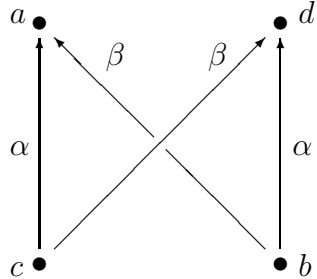

FIGURE 1: The graph $G$

The vertices labelled $a$ and $d$ in Figure 1 will be called the **top vertices** and the vertices labelled $b$ and $c$ will be called **bottom vertices**. For a given positive integer $n$, let $nG$ denote the graph comprised of $n$ disjoint copies of $G$. Let $\mathcal{M}$ denote a matching from the $2n$ top vertices to the $2n$ bottom vertices in $nG$. We consider the edges in $\mathcal{M}$ to be *downward* directed edges. Let $e$ denote a distinguished edge in $\mathcal{M}$.

**Definition 2.7** Given a data sequence $\langle n, \mathcal{M}, e \rangle$, a **restricted labelling** of the graph $nG$ is any labelling of the vertices in which each copy of $G \subseteq nG$ is labelled with the elements of a twisted $\alpha, \beta$–matrix (in the pattern shown in Figure 1), and where for each edge $f \in \mathcal{M} - \{e\}$ the head and tail of $f$ have the same label.

There is an apparent asymmetry between $\alpha$ and $\beta$ right now, since $\beta$–edges appear along the diagonal of each copy of $G$ and $\alpha$–edges do not. However this asymmetry is fictitious. A matrix $M$ belongs to $M(\alpha, \beta)$ if and only if the transpose of $M$ belongs to $M(\beta, \alpha)$, and therefore a matrix $M'$ belongs to $TM(\alpha, \beta)$ if and only if the matrix obtained from $M'$ by interchanging the bottom two elements belongs to $TM(\beta, \alpha)$. Therefore "top vertex", "bottom vertex" and "restricted labelling" have meanings which are unchanged if we switch the roles of $\alpha$ and $\beta$. (One can use this fact and the next lemma to give a new proof that the linear commutator is symmetric.)

**LEMMA 2.8** In $\mathbf{A}$ we have $[\alpha, \beta]_\ell = 0$ if and only if for each data sequence $\langle n, \mathcal{M}, e \rangle$ and each restricted labelling of the graph $nG$, it is the case that the head and tail of the distinguished edge are equal.

**Proof:** We have already seen that $[\alpha, \beta]_\ell = 0$ holds if and only if the following condition is satisfied.

> Given any (finite) sum of the form $\sum(a_i - b_i - c_i + d_i)$, where each summand comes from an $\alpha, \beta$–matrix, and where each element with a plus sign is matched with an identical element with a minus sign except that two elements $-u$ and $+v$ are left over, it is the case that $u = v$.

Assume that this condition is met and that we have a data sequence $\langle n, \mathcal{M}, e \rangle$ and a restricted labelling of $nG$. Take the sum of all labels of vertices in $nG$ with the top vertices given a plus sign and the bottom vertices given a minus sign. Then, since the labelling is restricted, we get a sum of the form $\sum(a_i - b_i - c_i + d_i)$, where each summand comes from an $\alpha, \beta$–matrix, and each element with a plus sign is matched with an identical element with a minus sign except that the labels of the head and tail of the distinguished vertex are left over. These left over labels must be equal elements (of opposite sign) by our criterion for $[\alpha, \beta]_\ell = 0$, and this implies that our criterion concerning restricted labellings is met.

Conversely, assume that our criterion for restricted labellings is met. Assume we are given a sum $\sum(a_i - b_i - c_i + d_i)$, where each summand comes from an $\alpha, \beta$–matrix, and where each element with a plus sign is matched with an identical element with a minus sign except that two elements $-u$ and $+v$ are left over. If there are $n$ summands, then we can use all of the elements $a_i, b_i, c_i$ and $d_i$ to label a copy of $nG$ in such a way that each copy of $G \subseteq nG$ is labelled with the elements of a twisted $\alpha, \beta$–matrix (using the same order for labels as is indicated in Figure 1). The matching of elements with a plus sign to elements with a minus sign determines a matching of all but one of the top vertices of $nG$ to all but one of the bottom vertices. The remaining top vertex is labelled $v$ and the remaining bottom vertex is labelled $u$. We complete the matching by taking $e$ to be the edge from this special top vertex to the special bottom vertex. With this choice our labelling of the vertices is a restricted labelling, so the head and tail of $e$ must have the same label. This implies that $u = v$. Therefore the criterion on restricted labellings implies our criterion for $[\alpha, \beta]_\ell = 0$. This completes the proof. □

## 3   A Sufficient Condition For $[\alpha, \beta]_s = [\alpha, \beta]_\ell$

In this section we shall analyze what it means for an algebra $\mathbf{A}$ to have congruences $\alpha$ and $\beta$ such that $[\alpha, \beta]_s = 0 < [\alpha, \beta]_\ell$. We shall find that whenever this happens, then for at least one choice of $\delta = \alpha, \beta$ or $\alpha \wedge \beta$ it is possible to define a congruence $\rho$ on $\mathbf{A} \times_\delta \mathbf{A}$ which has a strange property. The hypothesis that "there is no such $\rho$" corresponding to $\alpha$ and $\beta$ will therefore imply that $[\alpha, \beta]_s = 0 \Longrightarrow [\alpha, \beta]_\ell = 0$. We shall find in the next section that if $\mathbf{A}$ generates a variety which satisfies a nontrivial idempotent Mal'cev condition, then "there is no such $\rho$" for any $\alpha, \beta \in \text{Con}(\mathbf{A})$.

We begin by fixing an algebra $\mathbf{A}$ and congruences $\alpha, \beta \in \mathrm{Con}(\mathbf{A})$ such that $[\alpha, \beta]_s = 0 < [\alpha, \beta]_\ell$. We assume that $n$ is the least positive integer for which there is a data sequence $\langle n, \mathcal{M}, e \rangle$ and a restricted labelling of $nG$ which witnesses $[\alpha, \beta]_\ell \neq 0$ in the manner specified in Lemma 2.8. (This means that there is a restricted labelling of $nG$ where the head and tail of $e$ have different labels.) We leave it to the reader to verify that a failure for $n = 1$ is either a failure of $C(\alpha, \beta; 0)$ or of $C(\beta, \alpha; 0)$, and hence of $[\alpha, \beta]_s = 0$, so the least $n$ for which there is a failure is greater than one.

We fix a witness $\langle n, \mathcal{M}, e \rangle$ of $[\alpha, \beta]_\ell > 0$ and in this witness we shall call the copy of $G$ which contains the tail end of the distinguished edge $e$ the **critical square**. The (upward) $\alpha$–edges constitute a matching from the bottom vertices of $nG$ to the top vertices. These edges together with the (downward) edges in $\mathcal{M}$ determine a directed graph on the vertices of $nG$ in which every vertex has indegree one and outdegree one. It follows that this graph is a union of cycles, which we shall refer to as $\alpha$–**cycles**. Similarly, the $\beta$–**cycles** will be the cycles determined by $\mathcal{M}$ and the $\beta$–edges.

Notice that the edge preceding $e$ in the $\alpha$–cycle of $e$ is an $\alpha$–edge from the critical square. This implies that the $\alpha$–cycle containing $e$ contains at least one $\alpha$–edge from the critical square, although in fact it may contain both $\alpha$–edges from the critical square. A similar statement is true for $\beta$. We shall break our argument into cases according to which of the following conditions holds.

(I) The $\alpha$–cycle containing $e$ contains only one $\alpha$–edge from the critical square.

(II) The $\beta$–cycle containing $e$ contains only one $\beta$–edge from the critical square.

(III) The $\alpha$–cycle containing $e$ contains both $\alpha$–edges from the critical square and the $\beta$–cycle containing $e$ contains both $\beta$–edges from the critical square.

These are the only cases.

An edge $e' \in \mathcal{M} - \{e\}$ will play a central role in what follows. The choice of $e'$ is made differently in each of the three cases above. If we are in Case I, then the $\alpha$–cycle of $e$ contains only one $\alpha$–edge from the critical square, so only one bottom vertex in the critical square is on the $\alpha$–cycle of $e$. We denote by $e'$ the edge of $\mathcal{M}$ whose head is at the bottom vertex of the critical square which is *not* part of the $\alpha$–cycle of $e$. In particular, this means that the $\alpha$–cycle of $e$ is different from the $\alpha$–cycle of $e'$.

The way to choose $e'$ in Case I is depicted in Figure 2 where the distinguished edge $e$ is the edge from vertex $w$ to vertex $y$. The bottom vertex which does not belong to the $\alpha$–cycle of $e$ is vertex $x$. We have specified that $e'$ is the edge of $\mathcal{M}$ whose head is at $x$; it is therefore the edge from some top vertex (which we call $z$ in Figure 2) to the vertex $x$. Figure 2 does not show any $\beta$–edges.

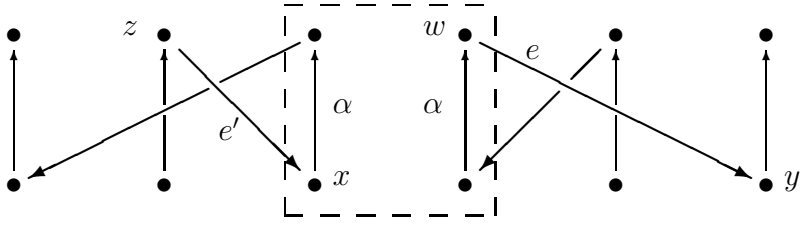

THE CRITICAL SQUARE

FIGURE 2: How to pick $e'$

If we are not in Case I but we are in Case II, then the choice of $e'$ is made as above after interchanging the roles of $\alpha$ and $\beta$. Thus $e'$ is the unique edge of $\mathcal{M}$ which has its head in the critical square, but which does not lie on the $\beta$–cycle of $e$.

Choosing $e'$ in Case III is a little more involved. This time the role that was played in Case I by the matching consisting of the $\alpha$–edges, and in Case II by the matching consisting of the $\beta$–edge will be played by a "mixed" matching selected as follows. We consider matchings from bottom vertices of $nG$ *upward* to top vertices where in each copy of $G$ we are free to choose *either* both $\alpha$–edges or both $\beta$–edges, but where these choices can be made independently in each copy of $G$. Among all such matchings, choose a matching $\mathcal{N}$ which maximizes the number of cycles in the graph whose edge set is $\mathcal{M} \cup \mathcal{N}$ and whose vertex set is the set of vertices of $nG$. Call the cycles formed by $\mathcal{M} \cup \mathcal{N}$ the $\zeta$–**cycles**. Observe that in no copy of $G$ is it the case that both $\alpha$–edges belong to the same $\zeta$–cycle, for if they did we could exchange these $\alpha$–edges from $\mathcal{N}$ for the $\beta$–edges in the same copy of $G$ and thereby increase the total number of cycles. Similarly, no two $\beta$–edges from the same copy of $G$ belong to the same $\zeta$–cycle. From this it follows that the $\zeta$–cycle of $e$ contains at most one bottom vertex from the critical square. We choose $e'$ to be the edge whose head is the bottom vertex of the critical square which is *not* in the $\zeta$–cycle of $e$.

**Definition 3.1** A **partially restricted labelling** of the graph $nG$ is any labelling of the vertices in which each copy of $G \subseteq nG$ is labelled with the elements of a twisted $\alpha, \beta$–matrix following the pattern in Figure 1, and where each edge in $\mathcal{M} - \{e, e'\}$ has the same label at its head and tail.

A restricted labelling is nothing more than a partially restricted labelling in which the head and tail of $e'$ have the same label.

Let $R$ denote the set of all quadruples $\langle (p, q), (r, s) \rangle$ for which $p, q, r$ and $s$ occur as the labels on vertices $w, x, y$ and $z$ in some partially restricted labelling of $nG$. This means the edge $e$ has labels $r$ and $p$ on its head and tail respectively and $e'$ has labels $q$ and $s$ on its head and tail respectively.

**LEMMA 3.2** *If $(a, b) R (c, d)$, then $a = b$ if and only if $c = d$.*

**Proof:** We shall prove that both of the possibilities $a = b$ and $c \neq d$, or else $a \neq b$ and $c = d$ contradict the minimality of $n$.

POSSIBILITY 1: $a = b$ AND $c \neq d$. Fix a partially restricted labelling of $nG$ satisfying these conditions. For now we assume that the head of $e$ and the tail of $e'$ do not belong to the critical square. Then our partially restricted labelling of $nG$ has the following form.

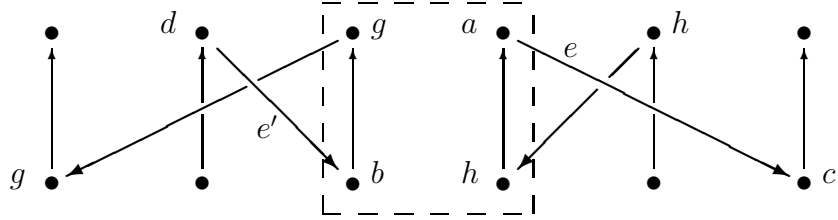

FIGURE 3

Observe that the multiple occurrences of labels $g$ and $h$ in this figure are forced since the head and tail of any edge in $\mathcal{M} - \{e, e'\}$ must receive the same label. Moreover we have $g = h$, since $a = b$, $[\alpha, \beta]_s = 0$ and the labels in the critical square are from a twisted $\alpha, \beta$–matrix in the pattern depicted in Figure 1. Now, by deleting the critical square we create a new data sequence $\langle n - 1, \widehat{\mathcal{M}}, \hat{e} \rangle$ which witnesses $[\alpha, \beta]_\ell > 0$. We let $\hat{e}$ be the edge whose tail is the old tail of $e'$ and whose head is the old head of $e$. We obtain $\widehat{\mathcal{M}}$ from $\mathcal{M}$ by deleting the four old edges which have a vertex in the critical square. These are $e, e'$, the edge with two $g$ labels, and the edge with two $h$ labels. Then we add $\hat{e}$ and a new edge from the top vertex with an $h$ label to the bottom vertex with a $g$ label.

The old partially restricted labelling yields a restricted labelling on our new graph by simply keeping all labels on the vertices that remain. Since $g = h$ this is indeed a restricted labelling. However, since the distinguished edge $\hat{e}$ has labels $c$ and $d$ on its head and tail, and $c \neq d$, we get from Lemma 2.8 that $\langle n - 1, \widehat{\mathcal{M}}, \hat{e} \rangle$ witnesses $[\alpha, \beta]_\ell > 0$ and this contradicts the minimality of $n$.

Other ways in which Possibility 1 might occur is if either $e$ or $e'$ has both its head and tail in the critical square. These subcases can be handled in the same way as above by slightly simpler arguments of the same form. We give the argument for the subcase where $e$ has both its head and tail in the critical square and $e'$ has its tail out of the critical square.

If the head of $e$ is in the critical square, then $e$ must be the reverse of either an $\alpha$–edge or a $\beta$–edge. Assuming the former, it is clear that the $\alpha$–cycle of $e$ is just $e$ followed by the reverse of $e$. (This means we are in Case I.) We have the following picture.

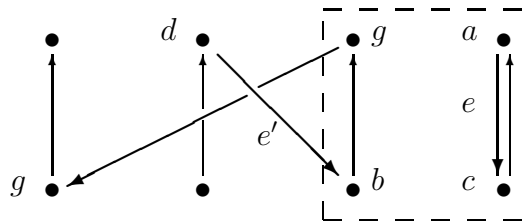

The same kind of argument as before shows that since $a = b$ and $[\alpha, \beta]_s = 0$, we have $c = g$. Thus, when we delete the critical square this time we choose $\hat{e}$ to be the new edge whose tail is the old tail of $e'$ and whose head is the bottom vertex which has label $g$. We obtain $\widehat{\mathcal{M}}$ from $\mathcal{M}$ by deleting the three old edges which have a vertex in the critical square, and then adding $\hat{e}$. The old partially restricted labelling of $nG$ yields a restricted labelling on the resulting graph where the labels on the distinguished edge $\hat{e}$ are again $c$ and $d$. Thus we get the same kind of contradiction as before.

The subcase where $e'$ has its tail in the critical square and the head of $e$ is outside can be handled by a symmetric argument. The case where both $e$ and $e'$ have their heads and tails in the critical square means that the four vertices of the critical square are labelled with equal elements on one edge of $G$ and unequal elements on the parallel edge. This is impossible since $[\alpha, \beta]_s = 0$.

POSSIBILITY 2: $a \neq b$ AND $c = d$. The argument here is the same as the argument for Possibility 1 with minor differences in detail. We argue only the main case, which is the one where the head of $e$ and the tail of $e'$ are outside the critical square. As in Possibility 1, we will delete the critical square to construct a data sequence $\langle n-1, \widehat{\mathcal{M}}, \hat{e} \rangle$ which witnesses $[\alpha, \beta]_\ell > 0$. Referring to Figure 3, observe that since $a \neq b$ and $[\alpha, \beta]_s = 0$ we must have that the labels $g$ and $h$ in the critical square are different. Therefore, if we let $\hat{e}$ be the edge from the top vertex labelled $h$ to the bottom vertex labelled $g$, and take $\widehat{\mathcal{M}}$ to be the matching obtained by deleting the four old edges which had vertex in the critical square and adding the edges $\hat{e}$ and a new edge from the tail of $e'$ to the head of $e$, then we have a data sequence $\langle n-1, \widehat{\mathcal{M}}, \hat{e} \rangle$ where the labels on the head and tail of $\hat{e}$ are $g$ and $h$ (which are unequal), but the labels on any edge in $\widehat{\mathcal{M}} - \{\hat{e}\}$ are the same at the head and tail.

We leave the other subcases under Possibility 2 to the reader. In all, the arguments show that if the statement of the lemma fails to hold then $n$ is not minimal. This contradicts our choice of $n$. $\square$

Let us set up notation which allows us to minimize arguing by cases. For the rest of this section we will let $\gamma$ and $\delta$ be equal to

- $\alpha$ and $\beta$ in Case I,
- $\beta$ and $\alpha$ in Case II,
- $\alpha \wedge \beta$ and $\alpha \wedge \beta$ in Case III.

This will not completely eliminate the need to consider cases separately and so the symbols $\alpha$ and $\beta$ will continue to be used. One thing worth noting at this point is that in all three cases $C(\delta, \gamma; 0)$ holds as a consequence of $[\alpha, \beta]_s = 0$.

**LEMMA 3.3** *The relation $R$ is a reflexive, compatible, binary relation on $\mathbf{A} \times_\delta \mathbf{A}$ for which there exist $(a,b) \, R \, (c,d)$ such that*

$$a \equiv_\gamma c \ \& \ b = d \ \& \ a \neq c.$$

**Proof:** To see that $R$ is a relation on $\mathbf{A} \times_\delta \mathbf{A}$ we must show that if $\langle (p,q), (r,s) \rangle \in R$, then $(p,q), (r,s) \in \delta$. Let us first argue this assuming that we are in Case I. The relevant edges of $nG$ and $\mathcal{M}$ are shown in Figure 4. (Recall that downward directed edges are edges in $\mathcal{M}$, while upward directed edges are either $\alpha$–edges or $\beta$–edges, according to the pattern set forth in Figure 1.) Since we are in Case I, $\delta = \beta$ and $(q,p)$ is a $\beta$–edge. Thus $(p,q) \in \beta = \delta$. For $(r,s)$ we argue as follows: beginning at $r$, the $\beta$–cycle of $e$ is $(r, \ldots, s, q, p)$. From $r$ to $s$ we have that consecutive labels are alternately equal (when they determine an edge from $\mathcal{M}$) and $\beta$–related (when they determine $\beta$–edges). This implies that $(r,s) \in \beta = \delta$. Thus $R$ is a relation on $\mathbf{A} \times_\delta \mathbf{A}$ in Case I. The Case II argument is the same with $\alpha$ and $\beta$ interchanged.

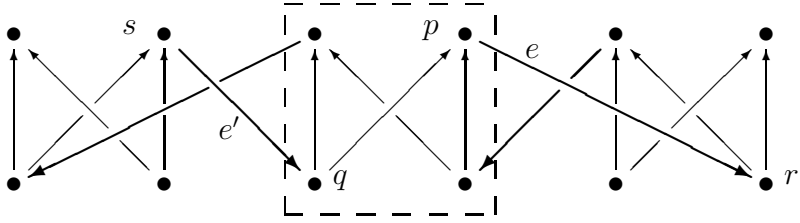

Figure 4

For Case III we must show that if $\langle (p,q), (r,s) \rangle \in R$, then $(p,q), (r,s) \in \delta = \alpha \wedge \beta$. The argument can again be followed in Figure 4, keeping in mind that now this figure depicts the situation when the $\mathcal{N}$–edges from the critical square happen to be $\alpha$–edges. Since $e'$ has its head in the critical square and the $\alpha$–cycle containing $e$ contains both $\alpha$–edges of the critical square, we get that $e'$ belongs to the $\alpha$–cycle of $e$. Similarly $e'$ belongs to the $\beta$–cycle of $e$. If $\langle (p,q), (r,s) \rangle \in R$, then there is some partially restricted labelling where the head and tail of $e$ are $r$ and $p$ respectively and the head and tail of $e'$ are $q$ and $s$ respectively. If we follow the $\alpha$–cycle of $e$ and $e'$ starting at $q$ we arrive at $p$ after traversing $\alpha$–edges going up and edges from $\mathcal{M} - \{e, e'\}$ going down, so $(p,q) \in \alpha$. The same argument starting at $r$ shows that $(r,s) \in \alpha$. We can apply the same argument along the $\beta$–cycle of $e$ and $e'$, and deduce that $(p,q), (r,s) \in \beta$. Thus $(p,q), (r,s) \in \gamma = \alpha \wedge \beta = \delta$. This finishes the proof that in all three cases $R$ is a binary relation on $\mathbf{A} \times_\delta \mathbf{A}$.

To see that $R$ is a compatible relation it suffices to show that $R$ is a subuniverse of $\mathbf{A}^4$. First observe that the set of all labellings of the $4n$ vertices of $nG$ with elements of $\mathbf{A}$ can be identified in a natural way with the algebra $\mathbf{A}^{4n}$. Since $TM(\alpha, \beta)$ is a subuniverse of $\mathbf{A}^4$, those labellings where each copy of $G$ is labelled with a twisted $\alpha, \beta$–matrix in the pattern of Figure 1 is a subuniverse which corresponds to $TM(\alpha, \beta)^n$. There is a subuniverse $S \subseteq A^{4n}$ consisting of all labellings of the vertices of $nG$ for which the head and tail of each edge in $\mathcal{M} - \{e, e'\}$ have the same label. Therefore the set of all partially restricted labellings, which is $S \cap TM(\alpha, \beta)^n$, is a subuniverse of $\mathbf{A}^{4n}$. The relation $R$ is just the projection of this subuniverse onto the four coordinates corresponding to the endpoints of $e$ and $e'$, so $R$ is a subuniverse of $\mathbf{A}^4$.

To show that $R$ is a reflexive relation on $\mathbf{A} \times_\delta \mathbf{A}$, we first assume that we are in Case I. In this case $\delta = \beta$ and therefore our task is to prove that for an arbitrary pair $(g,h) \in \beta$ it is the case that $\langle (g,h), (g,h) \rangle \in R$. To do this we need to show that there exists a partially restricted labelling of $nG$ such that the vertices $w, x, y$ and $z$ are labelled $g, h, g$ and $h$ respectively. We specify such a labelling by assigning $g$ to every vertex in the $\alpha$–cycle of $e$ and assigning $h$ to every other vertex. In this assignment the only labelled copies of $G$ that appear in $nG$ are copies with all labels equal to $g$, copies with all labels equal to $h$, and/or copies where two vertices along one $\alpha$–edge are labelled $g$ while the two vertices along the other $\alpha$–edge are labelled $h$. Each such labelling of a copy of $G$ is induced by a twisted $\alpha, \beta$–matrix. Moreover, since labels are constant along any $\alpha$–cycle, this is a (partially) restricted labelling. Thus, in Case I, we have established that $R$ is a reflexive relation on $\mathbf{A} \times_\delta \mathbf{A}$. The Case II argument is the same with $\alpha$ and $\beta$ interchanged.

To show that $R$ is reflexive in Case III, we must show that $\langle (g,h), (g,h) \rangle \in R$ whenever $(g,h) \in \delta = \alpha \wedge \beta$. We do this by labelling all vertices in the $\zeta$–cycle of $e$ with $g$ and all other vertices with $h$. Each copy of $G$ in $nG$ is labelled with all $h$'s or else with $g$'s on one upward directed edge and $h$'s on the parallel edge. Since $(g,h) \in \gamma = \alpha \wedge \beta$, each of these labellings of copies of $G$ is induced by a twisted $\alpha, \beta$–matrix. This way of assigning labels is constant on any $\zeta$–cycle, so it follows that every edge in $\mathcal{M}$ has the same label on its head and tail. Therefore this is a (partially) restricted labelling. Finally, since $e'$ is not on the $\zeta$–cycle of $e$ it follows that the head and tail of $e'$ are labelled $h$ while the vertices in the $\zeta$–cycle of $e$ are labelled $g$. Therefore we get $\langle (g,h), (g,h) \rangle \in R$ as desired.

What remains to prove is that there is a quadruple $\langle (a,b), (c,d) \rangle \in R$ such that $(a,c) \in \gamma - 0$ and $b = d$. The fact that we chose $n$ to witness $[\alpha, \beta]_\ell > 0$ implies that there is a restricted labelling of $nG$ such that all edges of $\mathcal{M} - \{e\}$ have the same label at the head and tail but that the labels on the head and tail of $e$ are different, say that the head and tail of $e$ are labelled $c$ and $a$ respectively. Our restricted labelling is simply a partially restricted labelling which has the property that the head and tail of $e'$ have the same label; call this label $b$. Taking $d = b$ shows that we have elements for which $\langle (a,b), (c,d) \rangle \in R$, $a \neq c$ and $b = d$. To finish the proof of this lemma we need only to show that $a \equiv_\gamma c$. Since $R$ is a relation on $\mathbf{A} \times_\delta \mathbf{A}$ we have $a \equiv_\delta b \equiv_\delta c$. In Case III we have $\gamma = \delta$, so there is nothing more to show here. In Case I, the $\alpha$–cycle of $e$ starting at $c$ may be written as $(c, \ldots, a)$ where the consecutive pairs of labels are alternately equal or $\alpha$–related. Since $\gamma = \alpha$ in this case, $a \equiv_\gamma c$. A similar argument works in Case II. The proof is finished. $\square$

We can combine the last two lemmas into a theorem about the situation where $[\alpha, \beta]_s = 0 < [\alpha, \beta]_\ell$. First we explain our notation for congruences in subalgebras of powers. The projection homomorphism from $\mathbf{A}^\kappa$ onto a sequence of coordinates $\sigma$ will be denoted $\pi_\sigma$. We will write $\eta_\sigma$ for the kernel of $\pi_\sigma$ and write $\theta_\sigma$ for $\pi_\sigma^{-1}(\theta)$ where $\theta$ is a congruence on $\mathbf{A}^\sigma$. The same symbols will be used for the restriction of these congruences to a subalgebra of $\mathbf{A}^\kappa$. The only exceptions to this rule are that 0 will denote the least congruence and 1 will denote the largest congruence.

**THEOREM 3.4** *Let $\mathbf{A}$ be an algebra with congruences $\alpha, \beta \in \mathrm{Con}(\mathbf{A})$. If $[\alpha, \beta]_s = 0 < [\alpha, \beta]_\ell$, then for at least one of the choices $(\gamma, \delta) = (\alpha, \beta), (\beta, \alpha)$ or $(\alpha \wedge \beta, \alpha \wedge \beta)$ there is a congruence $\rho \in \mathrm{Con}(\mathbf{A} \times_\delta \mathbf{A})$ such that*

(1) $\rho \leq \gamma_0 \wedge \gamma_1$,

(2) *the diagonal of $\mathbf{A} \times_\delta \mathbf{A}$ is a union of $\rho$–classes, and*

(3) $\rho \wedge \eta_1 \neq 0$.

**Proof:** We have shown in Lemmas 3.2 and 3.3 that if $[\alpha, \beta]_s = 0 < [\alpha, \beta]_\ell$, then it is possible to define a binary relation $R$ on $\mathbf{A} \times_\delta \mathbf{A}$ with the properties specified in those lemmas. The transitive closure, $\theta$, of $R \circ R^\cup$ is the congruence on $\mathbf{A} \times_\delta \mathbf{A}$ generated by $R$. We let $\rho = \theta \wedge \gamma_0 \wedge \gamma_1$. This definition of $\rho$ ensures that (1) holds.

The condition that $a = b \iff c = d$ whenever $(a, b)\ R\ (c, d)$ implies that the diagonal of $\mathbf{A} \times_\delta \mathbf{A}$ is a union of $\theta$–classes, and therefore of $\rho$–classes. This ensures that (2) holds.

The fact that there exist $(a, b)\ R\ (c, d)$ such that $a \equiv_\gamma c\ \&\ b = d\ \&\ a \neq c$ implies that $(a, b)$ and $(c, d)$ are distinct pairs for which

$$\langle(a, b), (c, d)\rangle \in \theta \wedge \gamma_0 \wedge \eta_1 = \rho \wedge \eta_1.$$

This ensures that (3) holds. $\square$

If $\mathbf{A}$ is any algebra and $\delta$ is any congruence on $\mathbf{A}$, then there is a largest congruence $\Delta$ on $\mathbf{A} \times_\delta \mathbf{A}$ such that the diagonal of $\mathbf{A} \times_\delta \mathbf{A}$ is a union of $\Delta$–classes. We denote this largest congruence $\Delta_\delta$. Notice that $\gamma_0 \wedge \gamma_1 \wedge \Delta_\delta$ is the largest congruence on $\mathbf{A} \times_\delta \mathbf{A}$ which satisfies the properties (1) and (2) attributed to $\rho$ in Theorem 3.4. If there is some $\rho$ as in Theorem 3.4 which satisfies (1), (2) and (3), then $\rho := \gamma_0 \wedge \gamma_1 \wedge \Delta_\delta$ is such a congruence. This yields the following result.

**THEOREM 3.5** *Let $\mathbf{A}$ be an algebra which has congruences $\alpha$ and $\beta$ such that $[\alpha, \beta]_s = 0$. Assume that for each choice of $(\gamma, \delta) = (\alpha, \beta),\ (\beta, \alpha)$ or $(\alpha \wedge \beta, \alpha \wedge \beta)$ we have*

$$\gamma_0 \wedge \eta_1 \wedge \Delta_\delta = 0$$

*in $\mathbf{Con}(\mathbf{A} \times_\delta \mathbf{A})$. Then $[\alpha, \beta]_\ell = 0$.*

An algebra is said to be **abelian** if it is abelian in the sense of the usual (TC) commutator, that is if $[1, 1] = 0$. This means the same thing as $[1, 1]_s = 0$. If $[1, 1]_\ell = 0$, then $\mathbf{A}$ is isomorphic to a subalgebra of a reduct of the affine algebra $\mathbf{A}^*/[1^*, 1^*]$ and so $\mathbf{A}$ is quasi–affine. The following corollary is the special case of the previous theorem where $\alpha = \beta = 1$.

**COROLLARY 3.6** *If $\mathbf{A}$ is abelian and $\Delta_1$ is a complement of the coordinate projection kernels in $\mathbf{Con}(\mathbf{A}^2)$, then $\mathbf{A}$ is quasi–affine.*

# 4 Imposing Mal'cev Conditions

The results in the previous section were local results in the sense that they were proved for individual algebras. In this section we derive global results from those local results by proving that there are nonobvious relationships between the affine, quasi–affine and abelian properties in varieties satisfying Mal'cev conditions. We then turn these results around to prove new results about Mal'cev conditions.

Our main result (Corollary 4.5) is that the symmetric commutator agrees with the linear commutator in any variety which satisfies a nontrivial idempotent Mal'cev condition. The most obvious consequence is the fact that in any variety satisfying a nontrivial idempotent Mal'cev condition the abelian algebras are quasi–affine. We further prove that: congruence neutrality is equivalent to congruence meet semidistributivity (Corollary 4.7), having a (weak) difference term is equivalent to a Mal'cev condition (Theorem 4.8 and the remarks that follow its proof), there are mild conditions on a variety under which one can conclude that abelian algebras are affine (Theorem 4.8, Theorem 4.10 and Corollary 4.11), and that every variety which satisfies a nontrivial lattice identity as a congruence equation has a weak difference term (Corollary 4.12).

We begin by defining what we mean by a Mal'cev condition. If $\mathcal{U}$ and $\mathcal{V}$ are varieties, then an **interpretation of $\mathcal{U}$ into $\mathcal{V}$** is a homomorphism from the clone of $\mathcal{U}$ to the clone of $\mathcal{V}$. If there is an interpretation of $\mathcal{U}$ into $\mathcal{V}$ we say that $\mathcal{U}$ is **interpretable** into $\mathcal{V}$ and we write $\mathcal{U} \leq \mathcal{V}$. If $\mathcal{U} \leq \mathcal{V}$ and $\mathcal{V} \leq \mathcal{U}$, then we write $\mathcal{U} \equiv \mathcal{V}$. The relation $\equiv$ is an equivalence relation on the class of varieties, and the $\equiv$–classes are called **interpretability classes**. Interpretability classes are ordered in a natural way by $\leq$, and under this order the collection of interpretability classes forms a complete lattice which we denote $\mathcal{L}$. It is known that the collection of all interpretability classes is a proper class, so our use of the phrase "complete lattice" for this collection only means that there are least upper bounds and greatest lower bounds with respect to $\leq$ for any set of interpretability classes.

Mal'cev conditions correspond to certain lattice filters in $\mathcal{L}$. To define these filters, let $\mathcal{C}$ denote any subset of interpretability classes. Now define a new (coarser) ordering on interpretability classes by saying that $u \leq_\mathcal{C} v$ whenever it is true that $w \leq u \implies w \leq v$ for all $w \in \mathcal{C}$. If $u \leq_\mathcal{C} v$ and $v \leq_\mathcal{C} u$, then write $u \equiv_\mathcal{C} v$ and say that $u$ and $v$ are **$\mathcal{C}$-indistinguishable**. We say that a filter in $\mathcal{L}$ is a **$\mathcal{C}$-filter** if it is a filter with regard to the $\leq_\mathcal{C}$–ordering. The proof of following easy lemma is left to the reader.

**LEMMA 4.1** *A collection $F$ of interpretability classes is a $\mathcal{C}$–filter if and only if*

(1) *it is a filter under $\leq$ and*

(2) *if $x \in F$ and $y$ is $\mathcal{C}$–indistinguishable from $x$, then $y \in F$.*

A $\mathcal{C}$–filter $F$ is said to be $(\mathcal{C}-)$**compact** if whenever $S \subseteq \mathcal{C}$ and $\bigvee S \in F$, then there is a finite subset $S_0 \subseteq S$ such that $\bigvee S_0 \in F$.

**Definition 4.2** A **Mal'cev filter** is a compact $\mathcal{C}$–filter where $\mathcal{C}$ is the collection of interpretability classes of finitely presentable varieties. A **Mal'cev condition** is an assertion of the form "the interpretability class of $\mathcal{V}$ belongs to $F$" where $F$ is a Mal'cev filter. The Mal'cev condition is **trivial** if $F = \mathcal{L}$ and **nontrivial** otherwise.

Our definition of 'Mal'cev condition' is formulated slightly differently than the usual definition found in, say, [5], but ours is an equivalent definition.

Our entire discussion concerning Mal'cev conditions could be carried out in the situation where $\mathcal{C}$ is the collection of interpretability classes of *idempotent*, finitely presentable varieties. In this case we obtain the definition of an **idempotent Mal'cev condition** as an assertion of the form "the interpretability class of $\mathcal{V}$ belongs to $F$" where $F$ is a compact idempotent Mal'cev filter. In the rest of this section we shall be concerned with varieties which satisfy a nontrivial idempotent Mal'cev condition. These are precisely the varieties $\mathcal{V}$ for which there is a finitely presented idempotent variety $\mathcal{I}$ such that $\mathcal{I} \leq \mathcal{V}$ but $\mathcal{I} \not\leq \mathcal{U}$ for some variety $\mathcal{U}$.

To begin our work we need the following lemma which characterizes those varieties which satisfy a nontrivial idempotent Mal'cev condition. This lemma is a direct consequence of Corollary 5.3 of [22].

**LEMMA 4.3** *A variety satisfies a nontrivial idempotent Mal'cev condition if and only if there is an $n > 1$, an idempotent $n$–ary term $f$ of $\mathcal{V}$ and $n$ linear equations satisfied in $\mathcal{V}$:*

$$\begin{aligned} f(x_{11}, \ldots, x_{1n}) &= f(y_{11}, \ldots, y_{1n}) \\ &\vdots \\ f(x_{n1}, \ldots, x_{nn}) &= f(y_{n1}, \ldots, y_{nn}) \end{aligned}$$

*where $x_{ij}, y_{ij}$ are variables and $x_{ii} \neq y_{ii}$ for each $i$.* □

Observe that in the previous lemma it is possible to choose all variables $x_{ij}, y_{ij}$ from the set $\{x, y\}$. For if we have a term $f$ which satisfies the kind of equations listed in the lemma we can specialize the variables to $\{x, y\}$ and still have equations of the same form: in equation $i$ we set whichever variable is in position $x_{ii}$ (and all other occurences of this variable in equation $i$) to the variable $x$. We set all other variables equal to $y$. This specialization is a consequence of the original equation, so it holds in $\mathcal{V}$, and we still have $x_{ii} \neq y_{ii}$.

**LEMMA 4.4** *Assume that $\mathcal{V}$ satisfies a nontrivial idempotent Mal'cev condition and that $\mathbf{A} \in \mathcal{V}$. If $\gamma$ and $\delta$ are congruences on $\mathbf{A}$ for which $C(\delta, \gamma; 0)$ holds, then in $\mathbf{Con}(\mathbf{A} \times_\delta \mathbf{A})$ we have*

$$\gamma_0 \wedge \eta_1 \wedge \Delta_\delta = 0.$$

**Proof:** Let $\theta = \gamma_0 \wedge \eta_1 \wedge \Delta_\delta$ and choose an arbitrary pair $\langle (a, c), (b, c) \rangle \in \theta$. Let $f(x_1, \ldots, x_n)$ be an idempotent term with the properties listed in Lemma 4.3. We can assume that $x_{11} = x \neq y = y_{11}$ and that the only variables in the equation $f(x_{11}, \ldots, x_{1n}) = f(y_{11}, \ldots, y_{1n})$ are $x$ and $y$. Substitute $b$ for all occurences of $x$ and $c$ for all occurences of $y$. This yields $f(b, \bar{u}) = f(c, \bar{v})$ where all $u_i$ and $v_i$ are in $\{b, c\}$. Since $b \equiv_\delta c$ we get that the polynomial defined by $p((x, y)) = (f(x, \bar{u}), f(y, \bar{v}))$ is a unary polynomial of $\mathbf{A} \times_\delta \mathbf{A}$. The equation $f(x_{11}, \ldots, x_{1n}) = f(y_{11}, \ldots, y_{1n})$ implies that $p((b, c))$ lies on the diagonal of $\mathbf{A} \times_\delta \mathbf{A}$. The element $p((a, c))$ is $\theta$–related to $p((b, c))$, and each element of the diagonal is a singleton $\theta$–class, therefore $p((a, c)) = p((b, c))$. This has the consequence that $f(a, \bar{u}) = f(b, \bar{u})$ where each $u_i \in \{b, c\}$. Now, since $\langle (a, c), (b, c) \rangle \in \theta \leq \gamma_1$ we get that $(a, b) \in \gamma$. Since $(a, c)$ and $(b, c)$ are elements of our algebra we have $a \equiv_\delta c \equiv_\delta b$. Therefore, applying $C(\delta, \gamma; 0)$ to the equality $f(a, \bar{u}) = f(b, \bar{u})$, we deduce that $f(a, \bar{y}) = f(b, \bar{y})$ for any $\bar{y} \in \{a, b\}^{n-1}$.

The argument we just gave concerning $a$, $b$ and $f$ which showed that $f(a, \bar{y}) = f(b, \bar{y})$ is an argument which works in each of the $n$ variables of $f$ if we choose the correct equation from Lemma 4.3. That is,

$$f(y_1, \ldots, y_{i-1}, a, y_{i+1}, \ldots, y_n) = f(y_1, \ldots, y_{i-1}, b, y_{i+1}, \ldots, y_n)$$

for each $i$ and any choice of values for $y_1, \ldots, y_n \in \{a, b\}$. Therefore, using the fact that $f$ is idempotent, we have $a = f(a, a, \ldots, a) = f(b, a, \ldots, ) = \cdots = f(b, b, \ldots, b) = b$. This proves that $\theta = 0$ and finishes the proof of the lemma. □

Now we are in a position to prove the main result of this section.

**COROLLARY 4.5** *If $\mathcal{V}$ is a variety which satisfies a nontrivial idempotent Mal'cev condition, then $\mathcal{V} \models [\alpha, \beta]_s = [\alpha, \beta]_\ell$. In particular, abelian algebras in $\mathcal{V}$ are quasi–affine.*

**Proof:** By Lemma 2.5 it suffices to show that $[\alpha, \beta]_s = 0 \implies [\alpha, \beta]_\ell = 0$. If $[\alpha, \beta]_s = 0$ holds, then $C(\delta, \gamma; 0)$ holds for each choice of $(\gamma, \delta) = (\alpha, \beta), (\beta, \alpha)$ and $(\alpha \wedge \beta, \alpha \wedge \beta)$. Therefore the hypotheses of Lemma 4.4 are met. The conclusion of Lemma 4.4 can then be used with Theorem 3.5 to conclude that $[\alpha, \beta]_\ell = 0$. □

It is not always easy to recognize if a variety satisfies a nontrivial idempotent Mal'cev condition so we now translate this condition into an equivalent one. We will say that an equation in the symbols $\{\vee, \wedge, \circ\}$ is a **congruence equation**. We intend to interpret the variables of the equation as congruence relations, $\vee$ as join of congruences, $\wedge$ as intersection of relations and $\circ$ as composition of (binary) relations. We say that a variety $\mathcal{V}$ **satisfies a congruence equation** if the equation holds in all congruence lattices of members of $\mathcal{V}$. A congruence equation is **trivial** if it holds in the congruence lattice of any algebra and **nontrivial** otherwise. Any congruence equation $u = v$ is equivalent to the pair of inclusions $u \subseteq v$ and $v \subseteq u$. Conversely, since our list of symbols includes $\wedge$, the congruence inclusion $u \subseteq v$ is equivalent to the congruence equation $u = u \wedge v$.

**LEMMA 4.6** *A variety satisfies a nontrivial idempotent Mal'cev condition if and only if it satisfies a nontrivial congruence equation.*

**Proof:** One direction of this proof is a standard argument and can be found in [17] or [23]. This is the direction which asserts that satisfaction of a nontrivial congruence equation implies the satisfaction of a nontrivial idempotent Mal'cev condition. The other direction is new so we include the proof.

Assume that $\mathcal{V}$ satisfies a nontrivial idempotent Mal'cev condition. By Lemma 4.3 we may assume that for some $n > 1$ the variety $\mathcal{V}$ has an idempotent $n$–ary term $f$ and that $\mathcal{V}$ satisfies
$$\begin{aligned} f(x_{11}, \ldots, x_{1n}) &= f(y_{11}, \ldots, y_{1n}) \\ &\vdots \\ f(x_{n1}, \ldots, x_{nn}) &= f(y_{n1}, \ldots, y_{nn}) \end{aligned}$$
where $x_{ij}, y_{ij} \in \{x, y\}$ for all $i$ and $j$ and $x_{ii} = x$, $y_{ii} = y$ for each $i$. We will use these equations to determine a nontrivial congruence equation for $\mathcal{V}$, so first we need some notation concerning these equations. Let $N = \{1, \ldots, n\}$. Let $L_i$ be the set of all $k \in N$ for which $x_{ik} = x$ and let $L'_i$ be the set of all $k \in N$ for which $x_{ik} = y$. Thus, $L_i$ and $L'_i$ describe the partition of $N$ which corresponds to the partition of the variables $\{x_{i1}, \ldots, x_{in}\}$ from the left hand side of the $i$–th equation into $x$'s and $y$'s. Let $R_i$ and $R'_i$ describe the partition on the right hand side of the equation: $R_i$ is the set of all $k$ for which $y_{ik} = x$ and $R'_i$ is the set of all $k$ for which $y_{ik} = y$. Now we are prepared to write down a congruence equation involving the variables $\{\alpha_1, \ldots, \alpha_n, \beta_1, \ldots, \beta_n\}$. So that the equation fits onto one line, let $\gamma = \bigwedge_N (\alpha_i \vee \beta_i)$ and let

$$\theta_i = ((\bigvee_{L_i} \alpha_i) \vee (\bigvee_{L'_i} \beta_i)) \wedge ((\bigvee_{R_i} \alpha_i) \vee (\bigvee_{R'_i} \beta_i)).$$

We claim that $\mathcal{V}$ satisfies the following congruence equation (which we write as an inclusion):

$$\bigwedge_N (\alpha_i \circ \beta_i) \subseteq ((\bigvee_N \alpha_i) \wedge \bigwedge_N (\gamma \vee \theta_i)) \vee ((\bigvee_N \beta_i) \wedge \bigwedge_N (\gamma \vee \theta_i)).$$

To see that $\mathcal{V}$ satisfies this congruence inclusion, choose $\mathbf{A} \in \mathcal{V}$ and congruences $\alpha_i$ and $\beta_i$, $1 \leq i \leq n$. Choose any $(a, b)$ from the relation defined by the left hand side of the inclusion. Then since $(a, b) \in \bigwedge_N (\alpha_i \circ \beta_i)$ we get that there exist $u_i \in A$ such that $a \equiv u_i \pmod{\alpha_i}$ and $b \equiv u_i \pmod{\beta_i}$. The element $U := f(u_1, \ldots, u_n)$ will play a crucial part in the argument. We make the following claims about the relationship between $U$ and the elements $a$ and $b$.

(1) $a \equiv U \pmod{\bigvee_N \alpha_i}$.

(2) $b \equiv U \pmod{\bigvee_N \beta_i}$.

(3) $a \equiv U \equiv b \pmod{\bigwedge_N (\gamma \vee \theta_i)}$.

Claim (1) is proved by noting that

$$a = f(a, \ldots, a) \equiv f(u_1, \ldots, u_n) = U \pmod{\bigvee_N \alpha_i}.$$

Claim (2) is proved the same way. For Claim (3), take the $i$–th equation satisfied by $f$ and substitute $a$ in for each occurrence of $x$ and $b$ for each occurrence of $y$. (That is, substitute $a$ for $x_{ik}$ and $y_{ik}$ if $k \in L_i \cup R_i$ and substitute $b$ for $x_{ik}$ and $y_{ik}$ if $k \in L'_i \cup R'_i$.) Let $v_i$ be the value obtained (on each side of the equation) after this substitution is made. From the left hand side of the equation we get $U = f(u_1, \ldots, u_n) \equiv v_i \pmod{((\bigvee_{L_i} \alpha_i) \vee (\bigvee_{L'_i} \beta_i))}$. From the right hand side of the equation we get $U \equiv v_i \pmod{((\bigvee_{R_i} \alpha_i) \vee (\bigvee_{R'_i} \beta_i))}$, therefore we have $U \equiv v_i \pmod{\theta_i}$. Since $v_i$ is obtained by substituting $a$'s and $b$'s into the arguments of $f$, and since $(a, b) \in \gamma$, we get that $a = f(a, \ldots, a) \equiv v_i \equiv f(b, \ldots, b) = b \pmod{\gamma}$. Altogether we get that

$$U \equiv v_i \equiv a \equiv b \pmod{\gamma \vee \theta_i}.$$

Since this holds for each $i$ we have that Claim (3) holds.

We can put Claims (1)–(3) together as follows. From Claims (1) and (3) we get that $a \equiv U \pmod{(\bigvee_N \alpha_i) \wedge \bigwedge_N (\gamma \vee \theta_i)}$. From Claims (2) and (3) we get that $U \equiv b \pmod{(\bigvee_N \beta_i) \wedge \bigwedge_N (\gamma \vee \theta_i)}$. This implies that

$$a \equiv b \mod ((\bigvee_N \alpha_i) \wedge \bigwedge_N (\gamma \vee \theta_i)) \vee ((\bigvee_N \beta_i) \wedge \bigwedge_N (\gamma \vee \theta_i)),$$

proving that $\mathcal{V}$ satisfies the congruence inclusion.

Now we must show that there is some algebra whose congruence lattice does not satisfy the congruence inclusion we are considering. The algebra we choose will be a set with no operations; specifically it will be the set $\{a, b, u_1, \ldots, u_n\}$. We define the following congruences (or equivalence relations) on this algebra. For $1 \leq i \leq n$ let $\alpha_i$ be the congruence with one nontrivial class, $\{a, u_i\}$, and let $\beta_i$ be the congruence with one nontrivial class, $\{b, u_i\}$. Observe that $(a, b) \in \bigwedge_N(\alpha_i \circ \beta_i)$. It can be easily checked that with these choices for $\alpha_i$ and $\beta_i$ the congruence $\gamma$ (as defined in the second paragraph of the proof) is equal to the congruence with one nontrivial class, $\{a, b\}$. Finally, from the equations for $f$, one sees that $x_{ii} = x \neq y = y_{ii}$, and this means that in the definitions of $L_i, L'_i, R_i$ and $R'_i$ we have that $i \in L_i \cap R'_i$ and $i \notin L'_i \cup R_i$. From this we get that $u_i$ is not $\theta_i$–related to either $a$ or to $b$. Consequently (since the only nontrivial $\gamma$–class is $\{a, b\}$) we get that $u_i$ is not $(\gamma \vee \theta_i)$–related to either $a$ or $b$. Thus none of the $u_i$'s are related to either $a$ or $b$ by the congruence $\bigwedge_N(\gamma \vee \theta_i)$. Since $a$ and $b$ are related to each other by this congruence (in fact by $\gamma$), we conclude that $\{a, b\}$ is a class of $\bigwedge_N(\gamma \vee \theta_i)$.

In this example the congruence $\bigvee_N \alpha_i$ has precisely two classes, which are $\{a, u_1, \ldots, u_n\}$ and $\{b\}$. The congruence $\bigvee_N \beta_i$ has only the classes $\{a\}$ and $\{b, u_1, \ldots, u_n\}$. From the result of the previous paragraph, we get that both $\{a\}$ and $\{b\}$ are singleton classes of $(\bigvee_N \alpha_i) \wedge \bigwedge_N(\gamma \vee \theta_i)$ and of $(\bigvee_N \beta_i) \wedge \bigwedge_N(\gamma \vee \theta_i)$. Therefore, they are singleton classes of $((\bigvee_N \alpha_i) \wedge \bigwedge_N(\gamma \vee \theta_i)) \vee ((\bigvee_N \beta_i) \wedge \bigwedge_N(\gamma \vee \theta_i))$. This proves that $(a, b)$ is in the left hand side of the congruence inclusion but not in the right hand side. $\square$

A variety is said to be **congruence neutral** if it satisfies the commutator congruence equation $[\alpha, \beta] = \alpha \wedge \beta$. It is **congruence meet semidistributive** if it satisfies the congruence implication $\alpha \wedge \beta = \alpha \wedge \gamma \implies \alpha \wedge \beta = \alpha \wedge (\beta \vee \gamma)$. To set up notation for the next result, assume that $\alpha, \beta$ and $\gamma$ are congruences. Define $\beta_0 = \gamma_0 = 0$, $\beta_{n+1} = \beta \vee (\alpha \wedge \gamma_n)$ and $\gamma_{n+1} = \gamma \vee (\alpha \wedge \beta_n)$

**COROLLARY 4.7** *Let $\mathcal{V}$ be a variety. The following conditions are equivalent.*

(1) $\mathcal{V}$ *is congruence neutral.*

(2) $\mathcal{V}$ *is congruence meet semidistributive.*

(3) $\mathcal{V}$ *satisfies the congruence equation $\alpha \wedge (\beta \circ \gamma) \subseteq \beta_n$ for some $n$.*

**Proof:** First we prove that (1) $\Rightarrow$ (2). Assume that $\mathcal{V}$ fails to be congruence meet semidistributive. Then we can find an algebra $\mathbf{A} \in \mathcal{V}$ which has congruences $\alpha, \beta$ and $\gamma$ such that
$$\delta := \alpha \wedge \beta = \alpha \wedge \gamma < \alpha \wedge (\beta \vee \gamma) =: \mu.$$
The part of this displayed line to the left of the "<" implies that $C(\beta, \alpha; \delta)$ and $C(\gamma, \alpha; \delta)$ hold, so $C(\beta \vee \gamma, \alpha; \delta)$ holds. Consequently $C(\mu, \mu; \delta)$ holds. This implies that $[\mu, \mu] \leq \delta < \mu$, which is a failure of congruence neutrality. Therefore (1) implies (2).

The equivalence of (2) and (3) is proved in [1]. So, to finish the proof it will suffice to prove that (2) implies (1). Suppose that $\mathcal{V}$ is congruence meet semidistributive. To prove that $\mathcal{V}$ is congruence neutral choose an arbitrary algebra $\mathbf{A} \in \mathcal{V}$ and congruences $\gamma = [\alpha, \beta] \leq \alpha \wedge \beta = \delta$. We must show that $\gamma = \delta$. Replacing **A** by a factor if necesary we may assume that $\gamma = 0$, so now we wish to prove that $\delta = 0$. Notice that $[\delta, \delta] \leq [\alpha, \beta] = 0$, so $\delta$ is a nonzero abelian congruence of **A**. This means that $C(\delta, \delta; 0)$, clearly.

By the equivalence of (2) and (3), and referring to (the easy direction of) Lemma 4.6, we get that $\mathcal{V}$ satisfies a nontrivial idempotent Mal'cev condition. From Lemma 4.4 we obtain that, since $C(\delta, \delta; 0)$ holds, then in $\mathbf{Con}(\mathbf{A} \times_\delta \mathbf{A})$ we have

$$(\delta_0 \wedge \eta_1) \wedge \Delta_\delta = 0.$$

By symmetry we also have $(\delta_1 \wedge \eta_0) \wedge \Delta_\delta = 0$, so by meet semidistributivity we also have

$$((\delta_0 \wedge \eta_1) \vee (\delta_1 \wedge \eta_0)) \wedge \Delta_\delta = 0.$$

But $(\delta_0 \wedge \eta_1) \vee (\delta_1 \wedge \eta_0) = \delta_0 \wedge \delta_1$, so we get $\delta_0 \wedge \delta_1 \wedge \Delta_\delta = 0$. Choose $(a, b) \in \delta$ arbitrarily. Since $C(\delta, \delta; 0)$ holds we have

$$\langle (a, a), (b, b) \rangle \in \delta_0 \wedge \delta_1 \wedge \Delta_\delta = 0,$$

and therefore $a = b$. This proves that $\delta = 0$ and concludes the proof. □

Corollary 4.7 (1) $\Leftrightarrow$ (2) is an extension (from the locally finite case to the general case) of a result that appears several times in different forms in [9]. (It appears in [9] for finite algebras in Corollary 5.20, for locally finite algebras in Lemma 7.7 (2) and for locally finite varieties in Theorem 9.10 (1)$\Leftrightarrow$(5).)

A **weak difference term** for a variety $\mathcal{V}$ is a ternary term $d(x, y, z)$ which has the property that whenever $\mathbf{A} \in \mathcal{V}$, $a, b, \in A$ and $\theta$ is a congruence on $\mathbf{A}$ such that $(a, b) \in \theta$, then

$$d(b, b, a) \; [\theta, \theta] \; a \; [\theta, \theta] \; d(a, b, b).$$

Clearly a weak difference term is Mal'cev on any block of an abelian congruence, so in particular any abelian algebra in a variety which has a weak difference term is affine.

**THEOREM 4.8** *Let $\mathcal{V}$ be a variety. The following conditions are equivalent.*

(1) $\mathcal{V}$ has a weak difference term.

(2) $\mathcal{V}$ satisfies the congruence equation $\alpha \wedge (\beta \circ \gamma) \subseteq (\alpha \wedge \beta_n) \circ \gamma \circ \beta \circ (\alpha \wedge \gamma_n)$ for some $n$.

(3) $\mathcal{V}$ satisfies an idempotent Mal'cev condition which is strong enough to imply that abelian algebras are affine.

**Proof:** Before starting on the proof, we elaborate on statement (3). What (3) says is that there is an idempotent Mal'cev filter $F$ containing the interpretability class of $\mathcal{V}$ which has the property that abelian algebras are affine in any variety whose interpretability class is in $F$.

The equivalence of (1) and (2) follows from Theorem 3.1 of [15] together with the results of this paper. Specifically, it is proved in [15] that (1) implies (2). The implication (2) $\Rightarrow$ (1) can be deduced from Theorem 3.1 of [15] provided we show that (2) implies that

$$[\alpha, \alpha] = 0 \implies [\alpha, \alpha]_\ell = 0$$

whenever $\alpha$ is a congruence on some algebra in $\mathcal{V}$. Now, (2) is a nontrivial congruence equation since it fails in the three–element set when we choose $\{\alpha, \beta, \gamma\}$ to be the three proper nonzero congruences. It follows from (the easy direction of) Lemma 4.6 that if (2) holds then $\mathcal{V}$ satisfies a nontrivial idempotent Mal'cev condition. Now clearly $[\alpha, \alpha] = 0$ is equivalent to $[\alpha, \alpha]_s = 0$; therefore Corollary 4.5 applies to show that the displayed implication is true whenever (2) holds. Thus we have (1) $\Leftrightarrow$ (2).

Any abelian algebra which has a weak difference term is affine, clearly. A variety which satisfies (2) satisfies an idempotent Mal'cev condition which by (1) is equivalent to having a weak difference term, and therefore any such variety has the property that its abelian algebras are affine. Conversely, assume that $\mathcal{V}$ satisfies an idempotent Mal'cev condition which implies that abelian algebras are affine. Then there is an idempotent Mal'cev filter $F$ which contains the interpretability class of $\mathcal{V}$ and which has the property that for every $\mathcal{U}$ whose interpretability class is in $F$ it is the case that the abelian algebras in $\mathcal{U}$ are affine. Take $\mathcal{U}$ to be the variety generated by the algebras which are idempotent reducts of algebras in $\mathcal{V}$. $\mathcal{U}$ satisfies any idempotent Mal'cev condition satisfied by $\mathcal{V}$, so $\mathcal{U} \in F$. If $\mathcal{A}$ is the subvariety of $\mathcal{U}$ generated by the abelian algebras, then $\mathcal{A}$ is an affine variety so there is a ternary term $t(x, y, z)$ such that $\mathcal{A} \models t(x, y, y) = x = t(y, y, x)$. We claim that $t$ is a weak difference term for $\mathcal{V}$. To see this, take $\mathbf{A} \in \mathcal{V}$ and $a, b \in A$. We must show that whenever $(a, b) \in \theta$ for some congruence $\theta$ then $t(a, b, b) \ [\theta, \theta] \ a \ [\theta, \theta] \ t(b, b, a)$. Factoring by $[\theta, \theta]$, one sees that it is enough to check this assertion when $\theta$ is an abelian congruence. But in this case the idempotent reduct $\mathbf{A}'$ of $\mathbf{A}$ is in $\mathcal{U}$ and the $\theta$–classes of $\mathbf{A}$ are abelian subalgebras of $\mathbf{A}'$, so they are in $\mathcal{A}$. Since $t(x, y, y) = x = t(y, y, x)$ in any algebra in $\mathcal{A}$, we get that $t(a, b, b) = a = t(b, b, a)$ in $\mathbf{A}$. This finishes the proof. □

A **difference term** for a variety $\mathcal{V}$ is a term $d(x, y, z)$ for which

$$d(b, b, a) = a \ [\theta, \theta] \ d(a, b, b)$$

holds whenever $(a, b) \in \theta$ and $\theta$ is a congruence of some member of $\mathcal{V}$. Using arguments from the previous proof and Theorem 4.1 of [15] one can prove that a variety has a difference term if and only if it satisfies the congruence equation $\alpha \wedge (\beta \circ \gamma) \subseteq \gamma \circ \beta \circ (\alpha \wedge \gamma_n)$ for some $n$. This fact is merely a digression. We proved Theorem 4.8 not because of a particular interest in (weak) difference terms, but simply because it gives us explicitly the weakest idempotent Mal'cev condition which implies that abelian algebras are affine: it is the Mal'cev condition corresponding to the congruence equation in Theorem 4.8 (2).

It is not always easy to recognize if a variety satisfies the congruence equation in Theorem 4.8 (2). Therefore we prepare to prove a result which is formally weaker than Theorem 4.8, yet which seems easier to apply. We prove that abelian algebras are affine in any variety satisfying an idempotent Mal'cev condition which fails to hold in the variety of semilattices. The following result is a direct consequence of Lemma 9.5 of [9].

**LEMMA 4.9** *A variety satisfies an idempotent Mal'cev condition which fails in the variety of semilattices if and only if there is an $n > 1$, an idempotent $n$–ary term $f$ of $\mathcal{V}$, and for every nonempty subset $K \subseteq \{1, \ldots, n\}$, there is an equation $f(x_{i_1}, \ldots, x_{i_n}) = f(y_{i_1}, \ldots, y_{i_n})$ satisfied in $\mathcal{V}$ where $\{x_{i_j} \mid j \in K\} \neq \{y_{i_j} \mid j \in K\}$ and the $x_{i_j}$ and $y_{i_j}$ are variables.* □

**THEOREM 4.10** *Assume that* **A** *generates a variety satisfying an idempotent Mal'cev condition which fails in the variety of semilattices. If* **A** *is abelian, then* **A** *is affine.*

**Proof:** We begin by replacing **A** by its idempotent reduct. The resulting algebra satisfies the same idempotent Mal'cev conditions as **A**, and it will be abelian if **A** is. Furthermore, an algebra is affine if and only if its idempotent reduct is affine. Therefore, we may assume from now on that **A** is an idempotent algebra.

We may invoke Corollary 4.5 to conclude that **A** is quasi–affine. If the similarity type of **A** is $\tau$, let $\mathcal{V}$ denote the variety of all $\tau$–algebras and let $\mathcal{V}^*$ be the variety of $\tau^*$–algebras defined in Section 2. There is a natural injective $\mathcal{V}$–homomorphism from **A** into $F(\mathbf{B})$ where **B** is the affine algebra $\mathbf{A}^*/[1^*, 1^*]$. If we show that $F(\mathbf{B})$ is affine, then it will follow that **A** is affine since $F(\mathbf{B})$ has a subalgebra isomorphic to **A**. Therefore, we lose no generality in assuming that **A** is in the image of the functor $F$, which implies that **A** is a reduct of an idempotent affine algebra, i.e., of an affine module. Since the affine module structure of **A** is generated by the operations of **A** together with $p(x, y, z) = x - y + z$, we may assume that the coefficients of the operations of **A** generate the coefficient ring **R**.

A two–generated free algebra in $\mathcal{V}(\langle \mathbf{A}; p(x, y, z)\rangle)$ can be constructed as the subalgebra of $\langle \mathbf{A}; p(x, y, z)\rangle^{|A|^2}$ generated by the projections. It is isomorphic to the affine **R**–module with universe $R$ and with generators 0 and 1. The reduct **B** of this algebra to the operations of **A** is an abelian algebra in $\mathcal{V}(\mathbf{A})$. Furthermore, the subalgebra of **B** generated by the projections (that is, generated with the operations of $\mathcal{V}(\mathbf{A})$) is a two–generated free algebra in $\mathcal{V}(\mathbf{A})$. Since **A** is affine if and only if $\mathbf{F}_{\mathcal{V}(\mathbf{A})}(2)$ is and **B** is affine if and only if all its subalgebras are affine, it follows that **A** is affine if and only if **B** is affine. Therefore we may replace **A** by an algebra isomorphic to **B** and from now on assume that **A** is a reduct of **R** considered as an affine **R**–module.

Let $I = \{c \in R \mid x - cy + cz$ is an operation of $\mathbf{A}\}$. It is not difficult to prove that $I$ is an ideal of **R** (see [20] or the proof of Lemma 5.9 of [11]). Of course, **A** is affine if and only if $1 \in I$. We assume now that this is not so and argue to a contradiction. Our first goal is to reduce to the case where $I = 0$. Let $\mathbf{S} = \mathbf{R}/I$ and let $\mathbf{C} = \mathbf{A}/\theta$ where $\theta = \{(r, s) \in A^2 = R^2 \mid r \equiv_I s\}$. Note that **C** is just the reduct to the operations of **A** of **S** considered as an affine **S**–module. The assumption that $I \neq R$ implies that **C** is nontrivial. Let $J = \{c \in S \mid x - cy + cz$ is an operation of $\mathbf{C}\}$. Choose $c \in J$. Since $x - cy + cz$ is an operation of **C**, there is an operation $ex - gy + hz$ of **A** such that $e \equiv_I 1$ and $g \equiv_I h$, where $c = g/I$. But now $(1 - e) \in I$, so we get that $x - (1-e)y + (1-e)z$ is also an operation of **A**. From $ex - gy + hz$ and $x - (1-e)y + (1-e)z$ we can construct

$$(ex - gy + hz) - (1-e)z + (1-e)x = x - gy + (e + h - 1)z = x - gy + gz$$

which must be an operation of **A**. Therefore $g \in I$ and so $c = g/I = 0$ in **S**. Since $c \in J$ was arbitrary, we get $J = 0$. In particular this shows that **C** is not affine.

Now we are in a position to use Lemma 4.9. There is an $n$–ary term $f$ witnessing the conclusions of Lemma 4.9 for the variety $\mathcal{V}(\mathbf{A})$, and therefore for the subvariety $\mathcal{V}(\mathbf{C})$. By permuting the variables if necessary we may assume that the affine expression for $f(z_1, \ldots, z_n)$ is $r_1 z_1 + \cdots + r_n z_n$ where for some $l \leq n$ we have that $r_1, \ldots, r_l$ are nonzero and $r_{l+1} = \cdots = r_n = 0$. Here $l \geq 1$ because $r_1 + \cdots + r_n = 1$. Let $K = \{1, \ldots, l\}$. The lemma states that there is an equation $f(x_{i_1}, \ldots, x_{i_n}) = f(y_{i_1}, \ldots, y_{i_n})$ satisfied in $\mathcal{V}$ where the set of variables $\{x_{i_1}, \ldots, x_{i_l}\}$ is different from $\{y_{i_1}, \ldots, y_{i_l}\}$. This implies that **C** satisfies

$$r_1 x_{i_1} + \cdots + r_l x_{i_l} = r_1 y_{i_1} + \cdots + r_l y_{i_l} \tag{1}$$

and that there is at least one variable that occurs on the left (say) which does not occur on the right. Call any such variable a **black** variable and call all other variables **white** variables. By further reordering the variables of $f$, we may assume that the black variables are precisely $x_{i_1}, \ldots, x_{i_k}$ for some $k$ such that $1 \leq k \leq l$. Now set all black variables equal to $y$ and all white variables equal to $x$. On the left hand side we get $r_1 y + r_2 y + \cdots + r_k y + sx$ where $s = r_{k+1} + \cdots + r_l$. On the right hand side we have only white variables, so this substitution yields only $x$. This specialization of equation (1) implies that

$$r_1 y + r_2 y + \cdots + r_k y + sx = x,$$

so $r_1 + \cdots + r_k = 0$ and $s = 1$. Now setting the first black variable on the left side of (1) equal to $z$, the rest of the black variables to $y$ and all white variables on the left to $x$ produces an operation of **C** which has the form $x + ty + r_1 z$ where $t = r_2 + \cdots + r_k = -r_1$. But now that $x - r_1 y + r_1 z$ is an operation of **C** we must have $r_1 = 0$, which is false. This contradiction completes the proof. □

Curiously, the following corollary is easy to derive from Theorem 4.10 but it does not seem to be easy to derive from Theorem 4.8, which is a stronger theorem. In this corollary a **lattice identity** is an equation in the symbols $\{\vee, \wedge\}$, and a lattice identity is **trivial** if it holds in every lattice and **nontrivial** otherwise.

**COROLLARY 4.11** (See Problem 1 of [9]) *If $\mathcal{V}$ satisfies a nontrivial lattice identity as a congruence equation, then the abelian algebras in $\mathcal{V}$ are affine.*

**Proof:** Let $\varepsilon$ be a lattice identity that holds in the congruence lattices of all algebras in $\mathcal{V}$, but which fails to hold in some lattice. It is proved in [17] and [23] that the collection of interpretability classes of congruence–$\varepsilon$ varieties is an intersection of idempotent Mal'cev filters. This implies that if a variety $\mathcal{S}$ is not congruence–$\varepsilon$, then there is an idempotent Mal'cev condition which holds in all congruence–$\varepsilon$ varieties but fails to hold in $\mathcal{S}$.

It is proved in [4] that if $\mathcal{S}$ is the variety of semilattices, then $\mathcal{S}$ is not congruence–$\varepsilon$ for any nontrivial $\varepsilon$. Therefore any congruence–$\varepsilon$ variety satisfies an idempotent Mal'cev condition which fails to hold in the variety of semilattices. Now just apply Theorem 4.10. □

**COROLLARY 4.12** *If $\mathcal{V}$ satisfies a nontrivial lattice identity as a congruence equation, then $\mathcal{V}$ has a weak difference term.*

We would like to close this section with an example. First, we have shown that if $\mathcal{V}$ satisfies a nontrivial idempotent Mal'cev condition, then abelian algebras in $\mathcal{V}$ are quasi–affine. Furthermore, if $\mathcal{V}$ satisfies an idempotent Mal'cev condition which fails in the variety of semilattices, then abelian algebras in $\mathcal{V}$ are affine. The question we have not yet answered is whether there is a variety which satisfies a nontrivial idempotent Mal'cev condition where the abelian algebras are *not* affine. We present such an example now.

The example we give solves Problem 3.6 of [15] which asks: if $\mathcal{V}$ satisfies a nontrivial congruence equation in the symbols $\{\vee, \wedge, \circ\}$ must $\mathcal{V}$ have a weak difference term? The answer is no, for the example we give satisfies a nontrivial idempotent Mal'cev condition, and therefore by Lemma 4.6 satisfies a nontrivial congruence equation in the symbols $\{\vee, \wedge, \circ\}$. However abelian algebras are not affine in this variety, so it cannot have a weak difference term according to Theorem 4.8.

**Example 4.13** Let $\mathbf{A}$ denote the reduct of the one–dimensional vector space over the real numbers to the operations of the form

$$r_1 x_1 + r_2 x_2 + \cdots + r_n x_n, \quad r_i \in [0,1] \text{ and } \Sigma r_i = 1.$$

Then $\mathbf{A}$ is not affine since any polynomial operation of $\mathbf{A}$ may be represented as a vector space polynomial with positive coefficients. Hence $x - y + z$ is not a polynomial operation of $\mathbf{A}$. However, $\mathcal{V}(\mathbf{A})$ satisfies a nontrivial idempotent Mal'cev condition, since the operation $f(x,y) := \frac{1}{2}x + \frac{1}{2}y$ fulfills the conditions listed in Lemma 4.3 for $n = 2$.

## 5 Related Remarks

For each affine algebra $\mathbf{A}$ there is a unique congruence on $\mathbf{A}^2$ which has the diagonal as a class, and this congruence is a complement of the coordinate projection kernels. Therefore, if $\mathbf{B}$ is quasi–affine then there is at least one congruence on $\mathbf{B}^2$ which has the diagonal as a class and which complements the coordinate projection kernels. We have proved in Corollary 3.6 that if $\mathbf{C}$ is abelian and *every* congruence on $\mathbf{C}^2$ which has the diagonal as a class is a complement of the coordinate projection kernels, then $\mathbf{C}$ is quasi–affine. This leaves open the following question: suppose that $\mathbf{C}$ is abelian and that there exists at least one congruence on $\mathbf{C}^2$ which has the diagonal as a class and which complements the coordinate projection kernels. Must $\mathbf{C}$ be quasi–affine? This question can be answered negatively in a very strong way, as we now explain.

Let $\tau$ be a similarity type and, following [18], write $\mathcal{A}(\tau)$ for the quasivariety of abelian algebras of type $\tau$ and $\mathcal{Q}(\tau)$ for the quasivariety of quasi–affine algebras of type $\tau$. Of course, $\mathcal{Q}(\tau) \subseteq \mathcal{A}(\tau)$, and this inclusion is proper when $\tau$ contains an operation of arity greater than one. In fact it is proved in [18] that $\mathcal{Q}(\tau)$ is not finitely based relative to $\mathcal{A}(\tau)$ unless the operations of $\tau$ are all unary. Now let $\mathcal{D}(\tau)$ be the class of those algebras $\mathbf{A}$ of type $\tau$ where $\mathbf{A}^2$ has a congruence $\Delta$ which complements the coordinate projection kernels in $\mathbf{Con}(\mathbf{A}^2)$ and which has the diagonal as a class. $\mathcal{D}(\tau)$, which is clearly contained in $\mathcal{A}(\tau)$, can be shown to be a quasivariety by showing that it is axiomatizable by a subset of the axioms in [18] which define $\mathcal{Q}(\tau)$ relative to $\mathcal{A}(\tau)$. Thus we have $\mathcal{Q}(\tau) \subseteq \mathcal{D}(\tau) \subseteq \mathcal{A}(\tau)$. The proof in [18] that $\mathcal{Q}(\tau)$ is not finitely based relative to $\mathcal{A}(\tau)$ when $\tau$ contains an operation of arity greater than one is actually a proof that $\mathcal{D}(\tau)$ is not finitely based relative to $\mathcal{A}(\tau)$. The question we raised in the previous paragraph may be restated as: is $\mathcal{Q}(\tau) = \mathcal{D}(\tau)$? The answer is no, since the arguments in [18] can be refined to show that $\mathcal{Q}(\tau)$ is not finitely based relative to $\mathcal{D}(\tau)$ when $\tau$ contains an operation of arity greater than one.

The results of this paper contribute something to the understanding of minimal idempotent varieties. Let $\mathcal{V}$ be any idempotent variety. If $\Lambda : \mathcal{V} \longrightarrow SETS$ is a clone homomorphism, then $\Lambda$ must be surjective. The kernel of $\Lambda$ is a set of defining equations for a subvariety of $\mathcal{V}$ which is definitionally equivalent to the variety of sets. Therefore, the following conditions are equivalent:

(1) $\mathcal{V}$ has no subvariety equivalent to the variety of sets.

(2) $\mathcal{V} \not\leq SETS$.

(3) there is a finitely presented idempotent variety $\mathcal{V}'$ such that $\mathcal{V}' \leq \mathcal{V}$ and $\mathcal{V}' \not\leq SETS$ (from (2), using compactness).

(4) $\mathcal{V}$ satisfies a nontrivial idempotent Mal'cev condition.

Similarly, $\mathcal{V}$ has no subvariety equivalent to the variety of sets or semilattices if and only if $\mathcal{V}$ satisfies a nontrivial Mal'cev condition which fails in the variety of semilattices.

Now let $\mathcal{V}$ be a minimal idempotent variety which is not equivalent to the variety of sets. If $\mathcal{V}$ is not congruence meet semidistributive, then by Corollary 4.7 $\mathcal{V}$ is not congruence neutral, so some algebra $\mathbf{A} \in \mathcal{V}$ has a pair of congruences $\alpha$ and $\beta$ such that $[\alpha, \beta] < \alpha \wedge \beta$. Then the congruence $\gamma := (\alpha \wedge \beta)/[\alpha, \beta]$ is a nonzero abelian congruence of $\mathbf{A}/[\alpha, \beta]$. If $\mathbf{B}$ is a subalgebra of $\mathbf{A}/[\alpha, \beta]$ supported by a nontrivial $\gamma$–class, then $\mathbf{B}$ is a nontrivial abelian algebra in $\mathcal{V}$. From the minimality of $\mathcal{V}$ and the existence of a nontrivial abelian algebra $\mathbf{B} \in \mathcal{V}$ we deduce that $\mathcal{V}$ has no subvariety equivalent to the variety of semilattices. Thus, Theorem 4.10 implies that $\mathbf{B}$ is affine. By the minimality of $\mathcal{V}$ we get that $\mathcal{V} = \mathcal{V}(\mathbf{B})$ is affine. Since an idempotent affine variety is minimal if and only if it is equivalent to a variety of affine modules over a simple ring, we have proved the following.

**THEOREM 5.1** *Let $\mathcal{V}$ be a minimal idempotent variety. $\mathcal{V}$ is equivalent to the variety of sets, to a variety of affine modules over a simple ring, or $\mathcal{V}$ is congruence meet semidistributive.*

This slightly improves a result in [11]. It is an open question as to which minimal idempotent varieties are congruence meet semidistributive, but when $\mathcal{V}$ is locally finite it is known from [21] that $\mathcal{V}$ must be congruence distributive or equivalent to the variety of semilattices.

**Acknowledgement.** The authors would like to thank Paolo Lipparini for showing us how Corollary 4.7 and Theorem 4.8 follow from our Corollary 4.5 and results in [15].

DEPARTMENT OF MATHEMATICAL SCIENCES
UNIVERSITY OF ARKANSAS
FAYETTEVILLE, AR 72701, USA.

BOLYAI INSTITUTE
ARADI VÉRTANÚK TERE 1
H–6720 SZEGED, HUNGARY.